\documentclass[10pt, reqno]{amsart}

\usepackage{amsmath, amscd, amssymb}
\usepackage[frame,cmtip,arrow,matrix,line,graph,curve]{xy}
\usepackage{graphpap, color}
\usepackage[mathscr]{eucal}

\newtheorem{theorem}{Theorem}[section]
\newtheorem{lemma}[theorem]{Lemma}
\newtheorem{proposition}[theorem]{Proposition}
\newtheorem{corollary}[theorem]{Corollary}
\theoremstyle{definition}

\newtheorem{conjecture}[theorem]{Conjecture}

\theoremstyle{remark}
\newtheorem{remark}[theorem]{Remark}

\numberwithin{equation}{section}

\newcommand{\Supp}{{\rm Supp}}

\newcommand{\vac}{|0\rangle}

{\vskip-\lastskip\medskip
  \noindent
  {\em #1.}\enspace
  }%
{\qed\par\medskip
  }

\def\CC{{\mathbb C}}
\def\QQ{{\mathbb Q}}
\def\RR{{\mathbb R}}

\def\cO{{\mathcal O}}
\def\cE{{\mathcal E}}
\def\cT{{\mathcal T}}
\def\cC{{\mathcal C}}

\def\cN{{\mathcal N}}
\def\cY{{\mathcal Y}}
\def\cU{{\mathcal U}}

\def\fM{\mathfrak{M}}
\def\fa{\mathfrak{a}}

\def\mapright#1{\,\smash{\mathop{\lra}\limits^{#1}}\,}

\def\eps{{\epsilon}}

\def\dual{^{\vee}}

\def\sta{^\ast}
\def\virt{^{\mathrm{vir}}}
\def\upmo{^{-1}}
\def\sta{^{\ast}}
\def\pri{^{\prime}}
\def\lam{{\lambda}}
\def\llam{_{\lambda}}
\def\spt{\text{Spt}}

\def\lra{\longrightarrow}
\def\lsta{_{\ast}}

\def\and{\quad{\rm and}\quad}
\def\bl{\bigl(}
\def\br{\bigr)}

\def\mh{\!:\!}
\def\sub{\subset}
\def\lab#1{\label{#1}[{#1}]\  }

\def\lalp{_\alpha}

\let\lab=\label

\def\mgnd{{\fM}_{g}(\xn,d)}

\def\xn{X^{[n]}}
\def\un{^{[n]}}

\def\Ob{\cO b}

\def\tilnul{}
\def\un{^{[n]}}
\def\on{^{(n)}}
\def\lfM{_{\fM}}
\def\on{^{(n)}}
\def\be{{\mathbf e}}
\def\bc{{\mathbf c}}
\def\bq{{\mathbf q}}

\def\image{\text{Im}}

\def\blangle{\bigl\langle}
\def\brangle{\bigr\rangle}

\begin{document}
\title[Gromov-Witten invariants of Hilbert schemes]
{two point extremal Gromov-Witten invariants of Hilbert schemes of
points on surfaces}

\author[Jun Li]{Jun Li$^1$}
\address{Department of Mathematics, Stanford University ,
Stanford,USA} \email{jli@math.stanford.edu}
\thanks{${}^1$Partially supported by an NSF grant DMS-0601002}

\author[Wei-Ping Li]{Wei-Ping Li$^2$}
\address{Department of Mathematics, HKUST, Clear Water Bay, Kowloon,
Hong Kong} \email{mawpli@ust.hk}
\thanks{${}^2$Partially supported by the grant CERG601905}


\begin{abstract}
Given an algebraic surface $X$, the Hilbert scheme $X^{[n]}$ of
$n$-points on $X$ admits a contraction morphism to the $n$-fold
symmetric product $X^{(n)}$ with the extremal ray generated by a
class $\beta_n$ of a rational curve. We determine the two point
extremal GW-invariants of $X^{[n]}$ with respect to the class
$d\beta_n$ for a simply-connected projective surface $X$ and the
quantum first Chern class operator of the tautological  bundle on
$X^{[n]}$. The methods used are vertex algebraic description of
$H^*(X^{[n]})$, the localization technique applied to $X=\mathbb
P^2$, and a generalization of the reduction theorem of Kiem-J. Li to
the case of meromorphic $2$-forms.
\end{abstract}

\maketitle
\date{}


\section{Introduction}

The Hilbert scheme $X^{[n]}$ of $n$-points of an algebraic surface
$X$ is a crepant resolution of the $n$-fold symmetric product
$X^{(n)}$.  The extremal ray of the contraction map $\pi\colon
X^{[n]}\to X^{(n)}$ is generated by the class of a rational curve
$\beta_n$ in $X^{[n]}$. The $k$-point extremal Gromov-Witten
invariants on $X^{[n]}$ is defined by
\begin{eqnarray*}
\blangle A^1, \ldots, A^k\brangle_{0,k, d}
=\int_{[\fM_{0,k}(\xn,d)]^{vir}}ev^*(A^1\otimes \ldots\otimes A^k),\
A^i\in H^*(X^{[n]},\mathbb C).
\end{eqnarray*}
When $X$ is a simply-connected projective surface, the $1$-point
extremal Gromov-Witten invariants of $X^{[n]}$ are  computed in
\cite{L-Q}. The main goal of this paper is to compute the $2$-point
extremal Gromov-Witten invariants of $X^{[n]}$.

\smallskip
Besides its own interest, this research is motivated by
the following two reasons. The first comes from a conjecture of Y.
Ruan. Since $X^{(n)}=X^n/S_n$ is an orbifold, one may ask a McKay
correspondence type question relating the cohomolgy ring   of
$X^{[n]}$ with the orbifold cohomolgy ring of $X^{(n)}$.  The
orbifold cohomology ring $H^*_{CR}(X^n/S_n)$ was defined by Chen and
Ruan in \cite{C-R} (see also \cite{AGV}). Using the extremal
Gromov-Witten invariants, Ruan defined a ``restricted" quantum
cohomology ring $H_\pi^*(X^{[n]})$ of the Hilbert scheme $X^{[n]}$
and conjectured in \cite{Ruan} that $H^*_{CR}(X^n/S_n)$ is
isomorphic to $H_\pi^*(X^{[n]})$ as rings (see also \cite{B-G}). The
orbifold cohomology ring $H^*_{CR}(X^n/S_n)$ was computed by
Fantechi-G\"ottsche and Uribe independently  in \cite{F-G, Ur}. Thus
the complete determination of the cohomology ring of $X^{[n]}$
depends upon the computation of the extremal Gromov-Witten
invariants on $X^{[n]}$.

\smallskip
The other motivation is the vertex algebraic nature of the
cohomology ring of $X^{[n]}$. Grojnowski \cite{Gro} and Nakajima
\cite{Na1} discovered that the direct sum of the cohomology groups
of $X^{[n]}$ over all $n$ is a highest weight irreducible
representation of a Heisenberg algebra.  From the work of Lehn
\cite{Lehn}, Lehn-Sorger \cite{L-S}, Li-Qin-Wang \cite{LQW}, and
Costello-Grojnowski \cite{C-G} on the cohomology ring of $X^{[n]}$,
the ring structure of  $X^{[n]}$ can be understood via Chern
characters of the tautological bundle on $X^{[n]}$ twisted by
cohomology classes of $X$. In particular, the first Chern class of
the tautological bundle plays a fundamental role. The right
viewpoint is that the first Chern class should be regarded as an
operator on the cohomology group $X^{[n]}$ via the cup product. The
operator  is determined by Lehn in \cite{Lehn} expressed in terms of
Heisenberg operators.  In the same spirit,  to understand the
quantum cohomology ring $H^*_\pi(X^{[n]})$, one needs to   study the
quantum first Chern class operator, which comes from the first Chern
class acting on the cohomology group of $X^{[n]}$    via the quantum
product in $H^*_\pi(X^{[n]})$. This is equivalent to computing the
two-point extremal Gromov-Witten invariants. When $X$ is the complex
plan equipped with a torus action, the full equivariant quantum
cohomology ring of $X^{[n]}$ is determined in \cite{O-P}.
The method used there is via localization, and the
key is to compute the quantum first Chern class operator.  In this
paper, we shall determine the quantum first Chern class operator
when $X$ is projective and simply-connected. As a consequence, the
two point extremal Gromov-Witten invariants are determined.

\smallskip
The main technique used in this paper is a generalized version of
the reduction theorem of the virtual cycle of the moduli space
$\fM_{g,k}(\xn,d)$ of   stable maps to the Hilbert scheme $X^{[n]}$.
The reduction theorem was first observed by Lee and Parker in
symplectic geometry. The algebro-geometric treatment of it was done
by Kiem and the first author in \cite{K-L}. The original theorem
deals with projective surfaces $X$ with a holomorphic two-form. In
order to cover general surfaces, we extend the
reduction theorem to cover the case where only meromorphic sections
of $\Omega_X^2$ are used. We now briefly outline the argument used
in this paper.

\smallskip

The Hilbert scheme $X^{[n]}$ contains a rational curve
\begin{eqnarray*}
\{\xi_{x_0}+x_1+\ldots+x_{n-2}\in X^{[n]}\mid
\Supp\{\xi_{x_0}\}=x_0\},
\end{eqnarray*}
where $x_0, x_1, \ldots, x_{n-2}$ are fixed points on $X$. This
curve is a generator of the extremal ray of the contraction map
\begin{eqnarray*}
\pi\colon X^{[n]}\longrightarrow X^{(n)}.
\end{eqnarray*}
Let $\beta_n$ represent the class of this curve in $H_2(X^{[n]},
\mathbb Z)$.

\smallskip
A stable map $\varphi\in \fM_{g,k}(\xn,d)$ can be factorized through
the product of punctual Hilbert schemes
\begin{eqnarray*}
\varphi=(\varphi_1,\cdots,\varphi_l): C\lra \prod
X^{[n_i]}_{p_i}\sub X\un,
\end{eqnarray*}
where   $\varphi_*(C)=d\beta_n$ and $\varphi_i$ is a morphism from
$C$ to the punctual Hilbert scheme $X^{[n_i]}_{p_i}$ consisting of
closed subschemes of $X$ of  length $n_i$ supported at a fixed point
$p_i$.

\smallskip
Suppose  $X$ admits a  holomorphic $2$-form $\theta$. The reduction
theorem basically says that the virtual cycle
$[\fM_{g,k}(\xn,d)]^{vir}$ is a homology class of a much smaller
space consisting of stable maps $\varphi$ satisfying that for each
$i$, either $\varphi_i$ is a constant map or its support $p_i$ lies
in the vanishing locus of $\theta$.  To extend the reduction theorem
to an arbitrary projective surface, we have to consider meromorphic
two-forms on $X$. The main part of the section \S 4 is to prove the
reduction theorem for a general surface. The conclusion is similar
except we replace the vanishing locus of a holomorphic two-form by
the zero locus and pole locus of a meromorphic two-form.

\smallskip
Using the reduction theorem, we conclude that  the
two-point extremal Gromov-Witten invariants satisfy a  universal
formula with universal coefficients to be determined. To get the
explicit expressions for the universal coefficients, we study the
equivariant two-point extremal Gromov-Witten invariants on the
projective plane $\mathbb P^2$ equipped with a torus action. Using
the localization, the computation is reduced to that on the complex
plane $\mathbb C^2$ which was done by Okounkov and Pandharipande in
\cite{O-P}. By the divisor axiom in the Gromov-Witten theory, we
also get the formula for the quantum first Chern class operator
expressed in terms of Heisenberg operators.

\smallskip
In this paper, all the homology and cohomology classes on $X^{[n]}$
are expressed in terms of Nakajima's  basis. Let $\mathfrak
a_i(\alpha)$ be the Heisenberg operators on the direct sum of the
cohomology of Hilbert schemes $X^{[n]}$ over all $n$, where $\alpha$
is a cohomology class on $X$.  Let $M(q)$ be the quantum first Chern
class operator $c_1(\mathcal O^{[n]})\cup_\pi $ where  $\mathcal
O^{[n]}$ is the tautological bundle on $X^{[n]}$ and $\cup_\pi$ is
the extremal quantum cup product defined via extremal Gromov-Witten
invariants. We have the following theorem.

\begin{theorem}
As an operator, the quantum first Chern class operator $M(q)$ can be
expressed in terms of Heisenberg operators as
\begin{eqnarray*}M(q) & =&\sum_{k>0}\big(\frac{k}{2}
         \frac{(-q)^k+1}{(-q)^k-1}-\frac{1}{2}\frac{(-q)+1}{(-q)-1}\big)
         \mathfrak a_{-k}\mathfrak a_k(\tau_{2*}[K_X])\nonumber\\
     &&~~~~~~~~~~~~~~~~~~~~~~~~~~~~~~~~    -\frac{1}{2}\sum_{k,\ell>0}\big(\mathfrak a_{-k}
\mathfrak a_{-\ell}\mathfrak a_{k+\ell}+\mathfrak
a_{-k-\ell}\mathfrak a_k\mathfrak a_{\ell}\big)(\tau_{3*}[X]).
\end{eqnarray*}
\end{theorem}

Again by the divisor axiom in the Gromov-Witten theory, the
two-point extremal Gromov-Wittin invariants of $X^{[n]}$ is thus
completely determined. As a corollary, Ruan's conjecture that
\begin{eqnarray*}
H^*_{CR}(X^n/S_n)\cong H^*_\pi(X^{[n]})
\end{eqnarray*}
is also verified for two-point case.

\smallskip
The paper is organized as follows. In  Section two, we review some
basic facts and notations about cohomology $H^*(X^{[n]})$ of the
Hilbert scheme $X^{[n]}$, especially Nakajima's treatment of
$H^*(X^{[n]})$. We review the extremal Gromov-Witten invariants of
Hilbert schemes. In Section three,  we generalize the reduction
technique via holomorphic two-forms of \cite{K-L} to the case of
meromorphic sections. In Section four, we prove the universality of
two-point extremal Gromov-Witten invariants of $X^{[n]}$  based on
the reduction theorem. As the result, the mentioned
Gromov-invariants only depend on certain universal coefficients. In
the next two Sections, we study the operator $M(q)$ and determine
these universal coefficients  by working on the projective plane
with a torus action. We achieve this by using the
localization technique to reduce to the computation of equivariant
Gromov-Witten invariants on the complex plane, which was already
done in \cite{O-P}. In Section seven, we determine the explicit
formula for the quantum first Chern class operator and the two-point
extremal Gromov-Witten invariants of the Hilbert scheme. We also
verify Ruan's conjecture for two-point case.

\smallskip
\noindent{\bf Acknowledgment}:  The first author would like to thank
Fudan University for the visit in 2006 where part of the research
was carried out.  The second author would like to thank MSRI at Berkeley, Stanford University and Fudan University for the visits in 2006 where the main part of the research was carried out.

\section{Preliminary on Hilbert schemes}
\subsection{The cohomology groups $H^*(X^{[n]}, \mathbb C)$ of the Hilbert scheme.}

For a smooth, simply-connected projective surface $X$,  we let
$X^{[n]}$ be the Hilbert scheme of length-$n$ zero dimensional
subschemes of $X$. The Hilbert scheme  $X^{[n]}$ is smooth; its
Hilbert-Chow morphism $\pi\colon X^{[n]}\to X^{(n)}$ is a crepant
resolution of the symmetric product $X^{(n)}$.

It is known that  the cohomology groups of $\xn$ form an infinite
dimensional vector space
\begin{eqnarray*}
\mathbb
H^X=\bigoplus_{n=0}^{\infty}\bigoplus_{k=0}^{4n}H^k(X^{[n]},\mathbb
C)
\end{eqnarray*}
that is the highest weight irreducible representation of a
Heisenberg algebra
$$\{\mathfrak a_i(\alpha)\}_{i\in\mathbb Z,
\alpha\in H^*(X)}
$$
of which the operators $\mathfrak a_i(\alpha)$ satisfy the
Heisenberg commutation relation
\begin{eqnarray}\label{Heisenberg}
[\mathfrak a_i(\alpha), \mathfrak a_j(\beta)]=-i\delta_{i+j,
0}\int_X \alpha\cup\beta,
\end{eqnarray}
and the highest weight vector
\begin{eqnarray*}
\vac=1\in H^0(\text{pt}, \mathbb C).
\end{eqnarray*}

We also know the generators of the vector space $H^*(X^{[n]},\mathbb
C)$: it is generated by elements of the forms
\begin{eqnarray*}
A^{\nu}=A^\nu(\alpha_1, \ldots, \alpha_r)=\mathfrak
a_{-\nu_1}(\alpha_1)\ldots\mathfrak a_{-\nu_r}(\alpha_r)\vac
\end{eqnarray*}
where $\nu\colon \nu_1\ge \nu_2\ge \ldots\ge\nu_r$ is a partition of
$n$ and  $\alpha_1, \ldots, \alpha_r$ are cohomology classes of $X$.

The same results hold for the homology groups of $X^{[n]}$. The
infinite dimensional vector space
\begin{eqnarray*}
\mathbb
H_X=\bigoplus_{n=0}^{\infty}\bigoplus_{k=0}^{4n}H_k(X^{[n]},\mathbb
C)
\end{eqnarray*}
is the highest weight irreducible representation of a Heisenberg
algebra
$$\{\mathfrak a_i(e)\}_{i\in\mathbb Z\mathbb e, e\in H_*(X)}
$$
where the operators $\mathfrak a_i(e)$ satisfy the Heisenberg
commutation relation
\begin{eqnarray}\label{Heisenberg2}
[\mathfrak a_i(e), \mathfrak a_j(e^\prime)]=-i\delta_{i+j, 0}\int_X
PD^{-1}e\cup PD^{-1}e^\prime,
\end{eqnarray}
and the highest weight vector
\begin{eqnarray*}
\vac=1\in H_0(\text{pt}, \mathbb C).
\end{eqnarray*}
The homology group $H_*(X^{[n]},\mathbb C)$ is also generated by
classes
\begin{eqnarray*}
A_\nu=A_{\nu}(e_1, \ldots, e_r)=\mathfrak
a_{-\nu_1}(e_1)\ldots\mathfrak a_{-\nu_r}(e_r)\vac.
\end{eqnarray*}
It is known that $PD^{-1}A_\nu=A^\nu$ if we take
$\alpha_i=PD^{-1}e_i$.

The homology class $A_\nu$ has a geometric description. If we
represent $e_i$ by geometric representatives $Z_i\sub X$, then
$A_\nu(e_1,\ldots, e_r)$ is a multiple of the closure of the set
\begin{eqnarray*}
\{\xi_1+\ldots+\xi_r\in X^{[n]}\, |\, \Supp(\xi_i)=x_i\in Z_i,
x_i\neq x_j \hbox{ for } i\neq j; \ \ell(\xi_i)=\nu_i\}
\end{eqnarray*}

One more useful remark: under the pairing $\blangle A^\lambda,
A^\mu\brangle=\int\limits_{X^{[n]}}A^{\lambda}\cup A^{\mu}$, the
operator $\mathfrak a_i(\alpha)$ is the adjoint operator of
$(-1)^i\mathfrak a_{-i}(\alpha)$.

A convention used through out the paper: for a subvariety $Y$ of a variety $X$, we use $[Y]$ to
represent the cohomology class in $H^*(X)$ dual to $Y$.

For the details of the results quoted here, the readers can consult
\cite{Na2}.

\subsection{Extremal Gromov-Witten invariants of Hilbert schemes}
\label{subsection-GW}

It is known that the Hilbert-Chow morphism $\pi\colon X^{[n]}\to
X^{(n)}$ contracts the curve class $\beta_n$ mentioned in the
introduction; the class $\beta_n$ also spans the extremal ray of the
morphism $\pi$ (see \cite{LQZ, V}).

In this paper, we shall investigate the moduli space of genus zero
stable morphisms to $\xn$ with  a $d$-multiple of the
extremal curve $\beta_n$ as the fundamental class;  we denote this moduli space by $\fM_{0, k}(X^{[n]}, d)$, where $k$ stands for the number of marked
points on the domains of the morphisms.
The $k$-point extremal Gromov-Witten invariants are defined via
\begin{eqnarray}\label{extremal-GW}
\blangle A^1, \ldots, A^k\brangle_{0, k, d}
=\int_{[\fM_{0,k}(\xn,d)]^{vir}}ev^*(A^1\otimes \ldots\otimes A^k),\
A^i\in H^*(X^{[n]},\mathbb C),
\end{eqnarray}
where $ev\colon \fM_{0,k}(\xn,d)\to
X^{[n]}\times\ldots\times X^{[n]}$ is the evaluation map at $k$
marked points.

Using the extremal Gromov-Witten invariants (\ref{extremal-GW}),
Ruan in \cite{Ruan} defined an extremal quantum cup-product
structure on $H^*(X^{[n]})$  as follows: by denoting
\begin{eqnarray*}
\blangle A^1, A^2, A^3\brangle_{qc}(q)=\sum_{d>0}\blangle A^1, A^2,
A^3\brangle_{0,3,  d}\cdot q^d
\end{eqnarray*}
and denoting
\begin{eqnarray*}
\blangle A^1, A^2, A^3\brangle_{qc}=\blangle A^1, A^2,
A^3\brangle_{qc}(q)|_{q=-1},
\end{eqnarray*}
he defined the quantum corrected cup product $A^1\cup_\pi A^2$ by
\begin{eqnarray*}
\blangle A^1\cup_\pi A^2, A^3\brangle=\blangle A^1\cup A^2,
A^3\brangle +\blangle A^1, A^2, A^3\brangle_{qc}.
\end{eqnarray*}
The cohomology group $H^*(X^{[n]})$ with so defined quantum product
$\cup_\pi$ is denoted by $H^*_\pi(X^{[n]})$.

 Chen and Ruan defined a ring structure on the
orbifold cohomology group  $H^*_{orb}(M/G)$ of the quotient of a
manifold $M$ by a finite group $G$. We denote this cohomology ring
by $H^*_{CR}(M/G)$. Applying this to the quotient $X\on=X^n/S_n$, we
obtain the Chen-Ruan orbifold cohomology ring $H\sta_{CR}(X\on)$.
The {\sl cohomological crepant resolution conjecture} formulated by
Ruan that relates the Chen-Ruan cohomology ring of an orbifold with
the quantum corrected cohomology ring of its crepant resolution, in
the case of Hilbert schemes, is of the following form:

\begin{conjecture}[Ruan]\label{Ruan}
$H^*_{CR}(X^{(n)})\cong H^*_\pi(X^{[n]})$ as rings.
\end{conjecture}

Since the ring $H^*_{CR}(X^{(n)})$ was computed in \cite{F-G, Ur}, to
verify this conjecture,  we need to derive an explicit form of the
quantum corrected cohomology ring of $H^*_\pi(X^{[n]})$, which will
be the task of the most of the remaining sections.

\section{The reduction Lemma}

Our approach relies on the reduction theorem of the virtual cycle of
the moduli space $\fM_{g,k}(\xn,d)$ observed by Lee-Parker in
symplectic geometry. The algebro-geometric treatment is the
localization technique worked out by Kiem and the first author in
\cite{K-L}. To begin with, we suppose the surface $X$ admits a
non-trivial holomorphic differential two-form
$\theta\in\Gamma(\Omega_X^2)$. Then following Beauville, $\theta$
induces a holomorphic two form $\theta\un$ of the Hilbert scheme
$\xn$, and by the result of \cite{K-L}, it defines a regular
cosection of the obstruction sheaf of $\fM=\mgnd$:
\begin{equation}
\lab{cosect} \eta: \Ob\lfM\lra \cO\lfM.
\end{equation}
Here $\Ob\lfM$ is the obstruction sheaf and $\cO\lfM$ is the
structure sheaf of $\fM$. We remark that marked points don't play
any role here. For simplicity, we only consider the case without
marked points.

This cosection reduces the virtual cycle of $\fM$ to a smaller
subset of it.

\begin{lemma}[\cite{K-L}] Let $\Lambda\sub\fM$ be the loci of points over which
$\eta$ fails to be surjective. Then the virtual cycle $[\fM]\virt\in
H\lsta(\Lambda)$.
\end{lemma}

To identify the vanishing loci of $\eta$, we recall the vanishing
criterion stated in \cite{K-L}: $\eta$ vanishes at $\varphi\in\fM$ if
the image of $\varphi\lsta: TC_{reg}\to T\xn$ lies in the null space
of
\begin{equation}\lab{theta}
\theta^{[n]}: T\xn\lra T\dual\xn.
\end{equation}
To pinpoint such $\varphi$, we notice that because of our choice of
the fundamental class of stable morphisms under investigation, the
composite of any $\varphi\in\fM$ with the Hilbert-Chow morphism
\begin{equation}\lab{HC}
\pi: \xn\lra X^{(n)}
\end{equation}
is a constant map. Let
$$\spt: \fM\lra X\on
$$
be the induced map. Then in case $\spt(\varphi)=\sum^l n_i p_i$, the
morphism $\varphi$ factors through the product of punctual Hilbert
schemes:
\begin{equation}\lab{decomp}
\varphi=(\varphi_1,\cdots,\varphi_l): C\lra \prod
X^{[n_i]}_{p_i}\sub X\un
\end{equation}
in which each $\varphi_i$ is a morphism from $C$ to the punctual
Hilbert scheme $X^{[n_i]}_{p_i}$. (Here, $X^{[m]}_p$ is the preimage
$\pi\upmo(mp)$ of $mp\in X^{(m)}$.) In the following, we call the
collection $\varphi=(\varphi_i)$ the standard decomposition of
$\varphi$ and call $p_i$ the support of $\varphi_i$. Note that the
collection $\{\varphi_i\}$ is canonical except the ordering; the
ordering depends on the ordering of points in $\spt(\varphi)$.

\begin{lemma}\lab{3.2}
Let $\Lambda_\theta\sub\fM$ be the set of those $\varphi\in\fM$
whose decompositions $\varphi=(\varphi_i)$ satisfying that for each
$i$ either $\varphi_i$ is a constant or its support
$p_i=\spt(\varphi_i)$ lies in the vanishing locus of $\theta$. Then
the locus where $\eta$ fails to be surjective is contained in
$\Lambda_\theta$.
\end{lemma}

\begin{proof}
Suppose $\varphi\in\fM$ lies in the vanishing loci of $\eta$ and
suppose $\varphi=(\varphi_i)$ is its decomposition. By the criterion
stated in \cite{K-L}, $\eta(\varphi)=0$ if and only if
$\varphi\lsta(TC_{reg})$ lies in the null space of $\theta\un$.
Because at a zero-dim subscheme $\xi$ that is a union of $l$
mutually disjoint zero-dim subscheme $\xi_i\in X^{[n_i]}$,
$$T_\xi \xn=\oplus_{i=1}^k T_{\xi_i}X^{[n_i]}
$$
and the form $\theta^{[n]}$ is a direct sum of $\theta^{[n_i]}$, the
image space $\varphi\lsta(TC_{reg})$ lies in the null space of
$\theta^{[n]}$ if and only if $\varphi_{i\ast}(TC_{reg})$ lies in the
null space of $\theta^{[n_i]}$ for all $i$. Now suppose $\xi_i$ is
supported at a single point $p_i$. By the work of Beauville, the
form $\theta^{[n_i]}$ is non-degenerate along $X^{[n_i]}_{p_i}$ if
$p_i\not\in \theta\upmo(0)$. Applying this to the support of
$\varphi_i$, we obtain the desired inclusion, thus proving the
Lemma.
\end{proof}

This reduction Lemma is sufficient for our application in case we
have a regular section $\theta\in H^0(X,K_X)$. For general surfaces,
we might not have such sections. Instead, we will work with
meromorphic sections of $K_X$ and show that such sections will
provide us the reduction lemma we need.

To this end, we let
\begin{equation}\lab{univ}
f: \cC\lra \xn\and \pi: \cC\to \fM
\end{equation}
be the universal family of $\fM$. As shown in \cite{L-T}, the
obstruction sheaf $\Ob\lfM$ of $\fM$ is a quotient sheaf of
$R^1\pi\lsta f\sta\cT_{\xn}$. We then pick a locally free sheaf
$\cE$ that surjects onto $R^1\pi\lsta f\sta T_{\xn}$, of which the
later surjects onto the obstruction sheaf $\Ob\lfM$. We let $E$ be
the vector bundle on $\fM$ whose sheaf of sections is $\cE$. Then
the construction of virtual cycle provides us a cone cycle $V\in
C\lsta E$ whose intersection with the zero section of $E$ gives rise
to the virtual cycle $[\fM]\virt$.

Next,  a meromorphic section $\theta$ of $K_X$,   viewed
as a meromorphic section of $\Omega_X^2$,  induces a meromorphic
section $\theta\un$ of $\Omega^2_{\xn}$, and hence a meromorphic
homomorphism
$$\eta: \cE\lra\cO_{\fM}.
$$
We let $D_0$ (resp. $D_\infty$) be the vanishing (resp. pole)
divisor of $\theta$.

Adopting the proof of the previous lemma, we immediately see that
the degeneracy loci of $\eta$:
$$Deg(\eta)=\{\varphi\in\fM\mid \text{either
$\eta$ is undefined or fails to be surjective at $\varphi$}\,\}
$$
is contained in the set of all $\varphi=(\varphi_i)$ such that
either for some $i$ the support of $\varphi_i$ is contained in
$D_0\cup D_\infty$ or for each $i$ the map $\varphi_i$ is a
constant.

It is the purpose of the remaining section to prove the reduction
Lemma, which says that the virtual cycle $[\fM]\virt$ lies in a much smaller
set than $Deg(\eta)$.

\begin{proposition}\lab{red}
Let $\Lambda_\theta\sub\fM$ be the subset that consists of those
$\varphi\in\fM$ whose decompositions $\varphi=(\varphi_i)$ have the
property that for each $i$ either $\varphi_i$ is constant or the
support of $\varphi_i$ lies in the union $D_0\cup D_\infty$. Then
the virtual cycle $[\fM]\virt$ is supported in $\Lambda_\theta$.
\end{proposition}

We first investigate the behavior of $\eta$ over where it is
undefined. For this, we introduce a partition of the moduli stack
$\fM$ based on the standard stratification of $X^{(n)}$. Recall that
the standard stratification of $X^{(n)}$ is indexed by the set of
all partitions of $n$, and that to  each partition
$\lambda=(\lambda_1\geq\ldots\geq\lambda_l)$ the stratum $X^\lam$ is
$$X^\lam=\{z=\sum_{i=1}^l \lam_ix_i\in X^{(n)} \mid x_1,\cdots,x_l\in X\
\text{are distinct}\}.
$$
Using preimages of the support map $\spt\mh \fM\to X^{(n)}$
mentioned in (\ref{HC}), we obtain a partition of $\fM$ indexed by
$\lambda$: $\fM_\lam=\spt\upmo(X^\lam)$, each endowed with the
reduced stack structure.

To proceed, we shall split the maps in $\fM_\lam$ into $l$
individual maps. To achieve this, we need an ordering of the points
occurring in the support of elements in $X^\lam$. We let
$$\psi\llam: X^{(\lam_1)}\times\cdots\times X^{(\lam_l)}
\lra X^{(n)}
$$
be the map that sends $(\xi_1,\cdots,\xi_l)$ to $\sum\xi_i$. Within
the domain of $\psi\llam$, we let $ B\llam$ be the open subset of
all $(\xi_1,\cdots,\xi_l)$ such that the support of $\xi_i$ are
mutually disjoint. Clearly, $\psi\llam( B\llam)=X^\lam$, and $\psi\llam\mh B\llam\to X^\lam
$ is \'etale.

Using $ B\llam$, we form $\cU\llam$ and the projection
$$\jmath\llam:\cU\llam=\fM\times_{X^{(n)}} B\llam\lra \fM;
$$
we let
$$\tilnul f_\lam: \tilnul\cC_\lam\lra X\un,\quad
\tilnul\pi\llam:\tilnul\cC\llam\lra \cU\llam
$$
be the pull back to $\cU\llam$ of the universal family $f$ via
$\jmath\llam$. Because elements in $\cU\llam$ are
$$(\varphi,(\xi_1,\ldots,\xi_l))\in \fM\times_{X^{(n)}} B\llam
$$
with support $\spt(\varphi)=\sum \xi_i$, as in (\ref{decomp})
$\varphi$ canonically splits into $l$ maps
$(\varphi_1,\cdots,\varphi_l)$ so that the support of $\varphi_i$ is
$\xi_i$. Obviously, this construction can be carried over to the
family $\tilnul f\llam$. In this way, we obtain $l$ morphisms
$$\tilnul f_{\lambda,i}: \tilnul \cC_{\lam}\lra X^{[\lambda_i]}
$$
such that over each closed point $(\varphi,(\xi_i))\in\cU\llam$ the
morphism $\tilnul f_{\lam,i}$ is the $\varphi_i$ alluded before.

Next we look at the obstruction sheaf $\Ob\lfM$. Recall that in
constructing the virtual cycle we have picked a locally free sheaf
$\cE$ surjects onto the sheaf $R^1\pi\lsta f\sta\cT_{\xn}$. As shown
in \cite{L-T}, we can pick $\cE$ so that over each $\cU\llam$, we
have direct sum decomposition
\begin{equation}
\lab{ds1}
\jmath\llam\sta\cE=\tilnul\cE_{\lam,1}\oplus\cdots\oplus\tilnul\cE_{\lam,l}
\end{equation}
and surjective homomorphisms
$$\tilnul\cE_{\lam,i}\lra R^1\tilnul\pi_{\lam\ast} \tilnul f_{\lam,i}\sta\cT_{X^{[\lam_i]}}
$$
that fits into the following commutative diagram
\begin{equation}
\lab{ds2}
\begin{CD}
\jmath\llam\sta\cE @= \tilnul\cE_{\lam,1}\oplus\cdots\oplus\tilnul\cE_{\lam,l}\\
@VVV @VVV\\
\jmath\llam\sta R^1\pi{\lam\ast} f\sta\cT_{\xn}@=
R^1\tilnul\pi_{\lam\ast} \tilnul
f_{\lam,1}\sta\cT_{X^{[\lam_1]}}\oplus\cdots\oplus R^1\pi_{\lam\ast}
f_{\lam,l}\sta\cT_{X^{[\lam_l]}}
\end{CD}
\end{equation}

We next look at the meromorphic homomorphism $\eta$. Following its construction,
$\eta$ is the composite
$$\cE\lra R^1\pi\lsta f\sta\cT_{\xn}\lra R^1\pi\lsta f\sta\Omega_{\xn}\lra
R^1\pi\lsta \omega_{\cC/\fM}\cong \cO_{\fM}
$$
in which the second arrow is induced by applying $\theta\un$, the
third arrow  by $f\sta$ and the last isomorphism by Serre's duality.
Similarly, replacing $\cE$ by $\tilnul\cE_{\lam,i}$ and replacing
$f\sta\cT_\xn$ by $\tilnul f_{\lam,i}\sta\cT_{X^{[\lam_i]}}$, we
obtain a  homomorphism
$$\tilnul\eta_{\lam,i}: \jmath\llam\sta\cE \mapright{\text{pr}}
\tilnul\cE_{\lam,i}\lra R^1\tilnul\pi_{\lam\ast} \tilnul
f_{\lam,i}\sta \cT_{X^{[\lam_i]}}\lra \cO_{\cU_\lam}.
$$
These individual (meromorphic) homomorphisms fit into the identity
$$\jmath\llam\sta\eta=
\tilnul\eta_{\lam,1}+\cdots+\tilnul\eta_{\lam,l},
$$
over where all make sense.

We now let $\Lambda_{\lam,i}\sub\cU\llam$ be those
$(\varphi,(\xi_i))$ such that either $\varphi_i$ are constant or the
supports $\spt(\varphi_i)$ of $\varphi_i$ satisfy
$$\spt(\varphi_i)\cap
(D_0\cup D_\infty)\ne\emptyset.
$$
Mimicking the proof in \cite{K-L}, we immediately see that each
$\tilnul\eta_{\lam,i}$ is surjective away from $\Lambda_{\lam,i}$.
Then applying results in \cite{K-L}, we obtain

\begin{lemma}
Away from $\Lambda_{\lam,i}$, the pull back cone $\jmath\llam\sta
N\sub\jmath\llam\sta E$ is contained in the kernel of
$\tilnul\eta_{\lam,i}$. Namely,
$$\jmath\sta N|_{\cU\llam-\Lambda_{\lam,i}}
\sub\ker\{\tilnul\eta_{\lam,i}:\jmath\llam\sta E\lra \cO_{\cU}\}.
$$
\end{lemma}

\begin{proof}[Proof of Proposition \ref{red}]
First we shall transform the problem from stacks to schemes. Let
$\cN\sub \cE$ be the cone over $\fM$ whose intersection with the
zero section gives the virtual cycle $[\fM]\virt$. We let $\cN\lalp$
be the irreducible components of $\cN$ with $c\lalp$ its
multiplicities. For each $\alpha$, we let $\cT\lalp\sub\fM$ be the
image stack of the projection $\cN\lalp\to\fM$; we pick a proper
variety $T\lalp$ and a morphism
$$\phi\lalp:T\lalp\lra \cT\lalp\sub\fM
$$
so that $\phi\lalp$ is generically finite. Since $\fM$ has a
projective coarse moduli space, such $T\lalp$ does exist.

We then let $E\lalp$ be the pull back vector bundle $\phi\lalp\sta
\cE$, and let $N\lalp\sub E\lalp$ be the subvariety so that under the
projection $E\lalp\to E|_{\cY\lalp}$ the variety $N\lalp$ maps
generically finitely onto $\cN\lalp$. Since both $\cN\lalp$ and
$T\lalp$ are irreducible, such $N\lalp$ is unique. Finally, we let
$s\lalp$ be the zero section of $E\lalp$; let $d\lalp$ be the degree
of the map $\phi\lalp$, which is identical to the degree of
$N\lalp\to \cN\lalp$. Then
$$[\fM]\virt=\sum\lalp c\lalp d\lalp\upmo\cdot
\phi_{\alpha\ast}s\lalp\sta[N\lalp].
$$

The reduction Proposition will follow from the following reduction
statement for each of the class $s\lalp\sta[N\lalp]$:
\begin{equation}\lab{rd3}
s\lalp\sta[N\lalp]\in H\lsta(\Lambda\lalp,\QQ), \quad\text{where}\
\Lambda\lalp=\phi\lalp\upmo(\Lambda_\theta).
\end{equation}

We now prove this statement for any fixed $\alpha$. We first cover
$T\lalp$ by open subsets $V_a$ so that each of its image
$\phi\lalp(V_a)$ is contained in the image of some $\cU\llam\to
\fM$. By choosing $V_a$ small enough, we can assume that $V_a\to
\fM$ lifts to $\phi_{a\lam}\mh V_a\to\cU\llam$. Using this, we can
pull back the direct sum decomposition (\ref{ds1}) and the
meromorphic homomorphisms $\eta_{\lam,i}$:
$$\begin{CD}
E\lalp|_{V_a}\cong \phi_{a\lam}\sta\cE_{\lam,1}\oplus\cdots\oplus
\phi_{a\lam}\sta\cE_{\lam,l}@>{\oplus
\phi_{a\lam}\sta\eta_{\lam,i}}>>
\cO_{V_a}\oplus\cdots\oplus\cO_{V_a}.
\end{CD}
$$
We denote $E_{a,i}=\phi_{a\lam}\sta \cE_{\lam,i}$ and
$\eta_{a,i}=\phi_{a\lam}\sta\eta_{\lam,i}$.

Like what we have done in \cite{K-L}, we shall pick (smooth) almost
splittings of the above homomorphisms. To control the behavior of
these splittings near $\phi_{a\lam}\upmo(\Lambda_{\lam,i})$, we pick
a small (analytic) neighborhood $\Lambda_{\lam,i}^\eps$ of
$\Lambda_{\lam,i}\sub \cU\llam$ so that it deformation retracts to
$\Lambda_{\lam,i}$.

Over $V_a$, we then pick a smooth section $\delta_{a,i}\in
C^\infty(V_a, E_{a,i})$ so that
\begin{enumerate}
\item $\delta_{a,i}\upmo(0)\sub \phi_{a\lam}\upmo(\Lambda_{\lam,i}^\eps)$;
\item the support of $\delta_{a,i}$, which is the (analytic) closure of
$\{\delta_{a,i}\ne 0\}$, is disjoint from
$\phi_{a\lam}\upmo(\Lambda_{\lam,i})$;
\item for $w\in V_a$ with $\delta_{a,i}(w)\ne 0$,
$\eta_{a,i}(\delta_{a,i}(w))$ is a positive real number.
\end{enumerate}
By first picking a smooth section $h$ of $E_{a,i}$ over
$V_a-\phi_{a\lam}\upmo(\Lambda_{\lam,i})$ so that $\eta_{a,i}\circ
h=1$ and then multiplying it by a cut off function, we obtain the
desired smooth section $\delta_{a,i}$.

Because of our choice of $\delta_{a,i}$, the section
$$\delta_a=\sum_{i=1}^l \delta_{a,i}\in C^\infty(V_a, E\lalp)
$$
has the property that
\begin{enumerate}
\item $\delta_a(w)=0$ if and only if $\delta_{a,i}(w)=0$ for all
$i$;
\item in case $\delta_a(w)\ne 0$, then the fiber of $N\lalp$ over $w$,
$N\lalp|_w$, lies in the kernel of $\eta_{a,i}(w)\mh E\lalp|_w\to
\CC$ for all $i$ of which $\delta_{a,i}(w)\ne 0$.
\end{enumerate}
Therefore, away from $\delta_a\upmo(0)$ the cone $N\lalp|_{V_a}$
lies in the kernel of $\delta_a\mh E\lalp|_{V_a}\to \CC_{V_a}$.

Our last step is to pick a partition of unity $g_a\mh V_a\to
\RR^{\geq 0}$ of the covering $\{V_a\}$; namely $\{g_a>0\}\Subset
V_a$ (its closure in $V_a$ is compact) and $\sum_a g_a\equiv 1$ on
$T\lalp$. The sum
$$\delta\lalp=\sum_a g_a\cdot\delta_a
$$
is then a smooth section of $E\lalp$. It is direct to check that
away from $\delta\lalp=0$, the cone $N\lalp$ is disjoint from
$\delta\lalp$. This proves that
$$s\lalp\sta[N\lalp]\in H\lsta(\delta\lalp\upmo(0)).
$$

It remains to pinpoint the set $\delta\lalp\upmo(0)$. This time,
because $\eta_{a,i}(\delta_{a,i}(w))\geq 0$ whenever it makes sense,
$\delta(w)\ne0$ if and only if for some $a$: $g_a(w)> 0$. Then
$\delta\lalp(w)=0$ implies $\delta_{a,i}(w)=0$ for all $i$, which
imply that $w\in \cap_i\phi_{a\lam}\upmo(\Lambda_{\lam,i}^\eps)$.
Therefore, $\delta\lalp\upmo(0)\sub \cup_a
\phi_{a\lam}\upmo(\Lambda_{\lam}^\eps)$. Finally, because of our
choice of $\Lambda_{\lam,i}^\eps$, we can retract $\cup_a
\phi_{a\lam}\upmo(\Lambda_{\lam}^\eps)$ to
$\Lambda\lalp=\phi\lalp\upmo(\Lambda_\theta)$ in $T\lalp$. This
proves that $s\lalp\sta[N\lalp]\in H\lsta(\Lambda\lalp)$. This
proves the reduction Proposition.
\end{proof}

\section{Universality of two point Gromov-Witten invariants}
\label{universality}

In this section, we will use the technique developed in the previous
section to prove a structure result on the two-point extremal
Gromov-Witten invariants of a general algebraic surface.

We consider the moduli space $\fM_{0, 2}(X^{[n]}, d)$ of two point genus zero
stable morphisms to $\xn$ of the $d$-multiple of
the extremal curve class $\beta_n$ as the  fundamental class. To determine the two point
Gromov-Witten invariants of this moduli space, we need to
investigate all possible
$$\blangle A_1, A_2\brangle_{0,2, d}^{\xn}\quad
\text{for}\ A_1, A_2\in H\sta(\xn).
$$
Using the reduction Proposition  of the previous section, we shall prove in
this section that,  with $A_i$ chosen among the Nakajima basis,  all but a
few such numbers vanish; and,  for those that don't,  their values only
depend on the intersection of $K_X$ with the relevant curve classes
involved.

To this end, we shall first recall the Nakajima basis of
$H\sta(\xn)$ and their geometric representatives. We let $\mu^1$,
$\mu^2$ and $\mu^3$ be three partitions of lengths $\ell(\mu^1)=r$,
$\ell(\mu^2)=s$ and $\ell(\mu^3)=t$; we write
\begin{eqnarray*}
&&\mu^1\colon \mu^1_1\ge \mu^1_2\ge\ldots\ge \mu^1_r,\quad
|\mu^1|=\mu^1_1+\mu^1_2+\ldots+\mu^1_r
\end{eqnarray*}
and write $\mu^2$ and $\mu^3$ accordingly. For the point class
$\bq\in H^4(X)$, curve classes $\bc_1,\cdots,\bc_s\in H^2(X)$ and
the fundamental class $[X]\in H^0(X)$, the triple
$\mu=(\mu^1,\mu^2,\mu^3)$ gives a cohomology class
\begin{eqnarray}\lab{Amu}
A^\mu_{\bc}=\mathfrak a_{-\mu^1_1}(\bq)\ldots\mathfrak
a_{-\mu^1_r}(\bq)\mathfrak a_{-\mu^2_1}(\bc_1)\ldots\mathfrak
a_{-\mu^2_s}(\bc_s)\mathfrak a_{-\mu^3_1}([X])\ldots\mathfrak
a_{-\mu^3_t}([X])\vac;
\end{eqnarray}
it is a class in $H^*(\xn)$ if $|\mu^1|+|\mu^2|+|\mu^3|=n$. By
going through all possible $\mu^i$ and classes $\bc_i$, the above
form a basis of the cohomology groups of $\xn$. To proceed, we keep
one such homology class $A^\mu_{\bc}$ as in (\ref{Amu}) and pick
three more partitions $\lambda^1$, $\lambda^2$ and $\lambda^3$ of
lengthes $a$, $b$ and $c$; pick curve classes $\be_1,\cdots,\be_b\in
H^2(X)$ and form
\begin{eqnarray}\lab{Ala}
A^\lambda_{\be}=\mathfrak a_{-\lambda^1_1}(\bq)\ldots\mathfrak
a_{-\lambda^1_a}(\bq)\mathfrak
a_{-\lambda^2_1}(\be_1)\ldots\mathfrak
a_{-\lambda^2_b}(\be_b)\mathfrak
a_{-\lambda^3_1}([X])\ldots\mathfrak a_{-\lambda^3_c}([X])\vac.
\end{eqnarray}
Again, it is a cohomology class of $\xn$ if
$|\lambda^1|+|\lambda^2|+|\lambda^3|=n$.

Our immediate task is to investigate the possibility of the
vanishing of $\blangle A^\lambda_{\be},A^\mu_{\bc}\brangle_{0,2,d}$.
For this, we need the  geometric representatives of the Poincar\'e
dual of $A^\lambda_{\be}$ and $A^\mu_{\bc}$. We pick points $q_1,
\ldots, q_r, p_1,\ldots, p_a$ in $X$; pick Riemann surfaces $C_i$
and $E_j$ that represent the Poincar\'e dual of the classes $\bc_i$
and $\be_j$, respectively. Without lose of generality, we can pick
these points and Riemann surfaces in general position that any
subcollection of them intersects transversally. Then the Poincar\'e
dual of $A^\mu_{\bc}$ is represented by
\begin{eqnarray*}
A_\mu^{\bc}=\mathfrak a_{-\mu^1_1}(q_1)\ldots\mathfrak
a_{-\mu^1_r}(q_r)\mathfrak a_{-\mu^2_1}(C_1)\ldots\mathfrak
a_{-\mu^2_s}(C_s)\mathfrak a_{-\mu^3_1}(X)\ldots\mathfrak
a_{-\mu^3_t}(X)\vac.
\end{eqnarray*}
$A_\mu^{\bc}$ is a multiple of the closure of the following subset
of $X^{[n]}$:

$$
\left\{\xi_1^1+\ldots+\xi^1_r+\xi^2_1+\ldots+\xi^2_s+\xi^3_1+\ldots+\xi^3_t
\,\Big|{{\xi^i_j\in X^{[\mu^i_j]},\,\Supp(\xi^3_j)=y_j\in
X}\atop{\Supp(\xi^2_j)=x_j\in C_j,\, \Supp(\xi^1_j)=q_j}}\right\}.
$$

Because the expected dimension of the moduli space $\fM_{0,
2}(X^{[n]}, d)$ is $2n-1$, $\blangle A^\lambda_{\be},
A^\mu_{\bc}\brangle_{0, 2, d}\neq 0$ is possible only if
$$\deg A^\lambda_{\be}+\deg A^\mu_{\bc}=\text{exp}.\dim \fM_{0, 2}(X^{[n]}, d)=2n-1.
$$
Then because the operators $\mathfrak a_{-k}(\bq)$, $\mathfrak
a_{-k}(\bc_i)$ and $\mathfrak a_{-k}([X])$ increase cohomology
degrees by $2k+2$, $2k$ and $2k-2$ respectively, the cohomology
degree of $A^\mu_{\bc}$ is $2(n+\ell(\mu^1)-\ell(\mu^3))$.
Therefore, the above identity forces
\begin{equation}
\big(\ell(\lambda^3)-\ell(\mu^1)\big)+
\big(\ell(\mu^3)-\ell(\lambda^1)\big)=1.\label{comp}
\end{equation}
Furthermore, by the reduction Proposition of the previous section,
$\blangle A^\lambda_{\be}, A^\mu_{\bc}\brangle_{0, 2, d}\neq 0$ only
if for the set $\Lambda_\theta$ defined there, with  $\pi\colon \xn\to X\on$,

\begin{equation}\lab{int}
\Lambda_\theta\cap \pi\upmo(A_\lambda^{\be})\cap
\pi\upmo(A_\mu^{\bc})\ne \emptyset,
\end{equation}

\begin{proposition}\label{reduction}
Suppose $d>0$ and $\blangle A^\lambda_{\be}, A^\mu_{\bc}\brangle_{0,
2, d}\neq 0$, then
\begin{equation*}
\ell(\lambda^3)=\ell(\mu^1)+\delta\quad\hbox{and}\quad
\ell(\mu^3)=\ell(\lambda^1)+1-\delta,\quad \text{for either} \ \
\delta=0 \ \text{or}\ 1.
\end{equation*}
In case $\delta=0$ holds, then  $\lambda^3=\mu^1$ as partitions; and
there exists an integer $\ell=\mu^3_i=\lambda^2_j$ for some integers
$i$ and $j$ such that the partition $\lambda^1$ is obtained from
$\mu^3$ with $\ell$ deleted, and the partition $\mu^2$ is obtained from
$\lambda^2$ with $\ell$ deleted.
\end{proposition}

\begin{proof}
Since $\blangle A^\lambda_{\be}, A^\mu_{\bc}\brangle_{0, 2, d}\neq
0$, there is an $(f,\Sigma)$ in the intersection (\ref{int}). We let
its support be the zero-cycle
$$\spt(f)=m_1x_1+\cdots+m_kx_k,\quad x_1,\cdots,x_k\quad
\text{distinct}.
$$
Since $f\in\pi\upmo(A_\lambda^{\be})$, there are three maps $u_1\mh
[a]\to [k]$, $u_2\mh [b]\to[k]$ and $u_3\mh [c]\to [k]$, where $[k]$
is the set of integers $\{1,\cdots,k\}$, such that $x_{u_1(i)}=p_i$,
$x_{u_2(i)}\in E_i$, that the coproduct
\begin{equation} \lab{cop}
u_1\sqcup u_2\sqcup u_3: [a]\coprod [b]\coprod [c]\lra [k]
\end{equation}
is surjective, and that
\begin{equation}\lab{B1}
m_l=\sum_{i\in u_1\upmo(l)}\lambda^1_i+\sum_{i\in
u_2\upmo(l)}\lambda^2_i+ \sum_{i\in u_3\upmo(l)}\lambda^3_i.
\end{equation}
For the same reason, since $f\in\pi\upmo(A_\mu^{\bc})$, we have maps
$v_1$, $v_2$ and $v_3$ from $[s]$, $[t]$ and $[r]$ to $[k]$,
respectively, such that $x_{v_1(i)}=q_i$, $x_{v_2(i)}\in C_i$, that
the coproduct $v_1\sqcup v_2\sqcup v_3$ is surjective, and that
\begin{equation}\lab{B2}
m_l=\sum_{i\in v_1\upmo(l)}\mu^1_i+\sum_{i\in v_2\upmo(l)}\mu^2_i+
\sum_{i\in v_3\upmo(l)}\mu^3_i.
\end{equation}

We first show that $\ell(\lambda^3)\geq\ell(\mu^1)$. Indeed, since
$q_1,\cdots,q_r$ are distinct, $v_1\mh [r]\to\image(v_1)$ is an
isomorphism. But then because of our general position requirement on
the points $p_i$ and $q_j$'s and of the Riemann surfaces $C_i$ and
$E_i$'s, $q_1, \ldots, q_r$ do not lie on $E_1,\ldots, E_b$,
$\image(v_1)\cap\bl \image(u_1)\cup\image(u_2)\br=\emptyset$. Hence
$\image(v_1)\sub\image(u_3)$ since $u_1\sqcup u_2\sqcup u_3$ is
surjective. This proves
$$\ell(\mu^1)=\# \image(v_1)\leq \#\image(u_3)\leq\ell(\lambda^3).
$$
For the same reason, we have $\ell(\mu^3)\ge \ell(\lambda^1)$.
Combined with (\ref{comp}), we obtain the first conclusion of the
Proposition.

Now suppose $\delta=0$; namely, $\ell(\mu^1)=\ell(\lambda^3)$, then
  all the identities above hold. In particular, $u_3$ is an
isomorphism $[c]\cong \image(u_3)\cong \image(v_1)$.

We next show that $\image(v_3)=\image(u_1)\cup\{k_0\}$ for an
integer $k_0\in [k]-\image(u_1)$. First, following the same reason
as before, we have $\image(u_1)\sub\image(v_3)$. Because
$\#\image(u_1)=\ell(\lambda^1)$ and because
$\ell(\lambda^1)+1=\ell(\mu^3)$, either $\image(v_3)=\image(u_1)$ or
$\image(v_3)=\image(u_1)\cup\{k_0\}$ for an $k_0\in
[k]-\image(u_1)$. We will show that only the later can happen.

For this, we decompose $f$ into $k$ individual morphisms $f_i\mh
\Sigma\to X^{[m_i]}_{x_i}$, where $X^{[m]}_x$ is the preimage of
$mx\in X^{(m)}$ under the Hilbert-Chow morphism $X^{[m]}\to
X^{(m)}$. Because $\pi(f(\Sigma))=\sum m_i x_i$, such decomposition
is possible. Since $d>0$, there is at least one $k_0$ so that
$f_{k_0}$ is non-constant. By the characterization of
$\Lambda_\theta$, this is possible only if $x_{k_0}\in D_0\cup
D_\infty$. Because of this, $x_{k_0}\ne p_i$'s, and thus
$k_0\not\in\image(u_1)$; also $k_0\not\in\image(u_3)$ because
$\image(u_3)=\image(v_1)$. Thus there is an $i_0\in [b]$ such that
$u_2(i_0)=k_0$. Thus $x_{k_0}\in E_{i_0}\cap(D_0\cup D_\infty)$,
which then exclude the possibility that $x_{k_0}\in C_j$'s. Hence
$k_0\in\image(v_3)$; and $\image(v_3)=\image(u_1)\cup\{k_0\}$. Note
that this also proves that $f_{k_0}$ is the only non-constant
$f_j$'s.

Now let $i\in [b]-\{i_0\}$ and consider $x_{u_2(i)}$. Because
$$\{x_l\mid l\in\image(v_1)\cup \image(v_3)\}\sub
\{q_1,\cdots,q_r\}\cup \{p_1,\cdots,p_a\}\cup \bl E_{i_0}\cap
(D_0\cup D_\infty)\br,
$$
$u_2(i)\not\in \image(v_1)\cup\image(v_3)$; hence there is a $j\in
[s]$ so that $u_2(i)=v_2(j)$, and consequently $x_{u_2(i)}\in
E_i\cap C_j$. Because $C_i$'s and $E_j$'s are in general positions,
once $x_{u_2(i)}\in E_i\cap C_j$, it does not lie in any other
$E_{i\pri}$'s and $C_{j\pri}$'s. In particular, $u_2\upmo\circ v_2$
defines an isomorphism $[s]\to [b]-\{i_0\}$.

Combined, we see that the map $u_1\sqcup u_2\sqcup u_3$ is
injective and  thus is an isomorphism. Similarly, $v_1\coprod v_2\coprod
v_3$ is also an isomorphism. Therefore, by (\ref{B1}) and
(\ref{B2}), we have
$$ \lambda^1_{u_1\upmo(i)}=\mu^3_{v_3\upmo(i)},\quad
\lambda^2_{u_2\upmo(i)}=\mu^2_{v_2\upmo(i)}\and
\lambda^3_{u_3\upmo(i)}=\mu^1_{v_1\upmo(i)}
$$
when $i\in\image(u_1)$, $i\in \image(u_2)-\{i_0\}$ and $i\in
\image(u_3)$, respectively.

Putting them together, we have proved,  in case $\delta=0$,  that $\lambda^3=\mu^1$ as
partitions, that there is an integer $\ell$ so
that $\lambda^2$ is $\mu^2$ with $\ell$ added, and that $\lambda^1$
is $\mu^3$ with $\ell$ deleted. Furthermore, the decomposition of $f$
has all but one component constant; the non-constant component is
the one associated to the part $\ell$.
\end{proof}

Let $A_{\be}^{\lambda-\lambda^2_j}$ be  the cohomology class on $X^{[n]}$ obtained from $A_{\be}^\lambda$  in (\ref {Ala}) with
$\mathfrak a_{-\lambda^2_j}( \be_j)$ deleted. Similarly, we can define $A_{\be}^{\lambda-\lambda^1_j}$, $A_{\be}^{\lambda-\lambda^2_i-\lambda^2_j}$, etc. For example,
\begin{eqnarray*}
A_{\be}^{\lambda-\lambda^2_1}=\mathfrak a_{-\lambda^1_1}(\bq)\ldots\mathfrak
a_{-\lambda^1_a}(\bq)\mathfrak
a_{-\lambda^2_2}(\be_2)\ldots\mathfrak
a_{-\lambda^2_b}(\be_b)\mathfrak
a_{-\lambda^3_1}([X])\ldots\mathfrak a_{-\lambda^3_c}([X])\vac.
\end{eqnarray*}

\begin{corollary} \lab{4.2}
Suppose $\lambda$ and $\mu$ fit into the
case $\delta=0$ in Proposition \ref{reduction}, then
$$
\blangle A^\lambda_{\be}, A^\mu_{\bc}\brangle_{0, 2, d} =\sum_{
\mu^3_i=\lambda_j^2} \blangle A^{\lambda-\lambda^2_j}_{\bc},
A^{\mu-\mu^3_i}_{\be}\brangle\cdot \blangle \mathfrak
a_{-\lambda^2_j}(\be_j)\vac , \mathfrak
a_{-\mu^3_i}([X])\vac\brangle_{0, 2, d}.
$$
\end{corollary}

\begin{proof} The proof is obvious and is omitted.
\end{proof}

Thus we only need to determine $\blangle \mathfrak
a_{-\ell}(\be)\vac , \mathfrak
a_{-\ell}([X])\vac)\brangle^{X^{[\ell]}}_{0, 2, d}$. For this we
have

\begin{lemma}\label{uni-constant}
There exists a universal function $c_{\cdot,\cdot}$ such that for
any positive integers $\ell$ and $d$, and homology class $\be\in
H^2(X)$,
\begin{eqnarray*}
\blangle \mathfrak a_{-\ell}(\be)\vac, \mathfrak
a_{-\ell}([X])\vac\brangle_{0, 2, d}=c_{\ell, d}\cdot( \be\cdot
c_1(K_X)).
\end{eqnarray*}
\end{lemma}

\begin{proof}
We first introduce the universal constant $c_{\ell,d}$. We let $U$
be a smooth analytic surface, let $\theta_+\in H^0(U, K_U)$ be an
analytic section vanishing along a smooth curve $C\sub U$, and let
$E\sub U$ be another smooth curve that intersects transversally with
$C$ at a single point $p\in U$. We can form the Hilbert scheme
$U^{[\ell]}$ as an analytic space and form the moduli of stable
morphisms $\fM_{0,2}(U^{[\ell]},d)$. By choosing $U$ as an analytic
open subset of a smooth algebraic surface, both $U^{[\ell]}$ and the
moduli space are analytic spaces of a projective scheme and of a
Deligne-Mumford stack; thus their existence are well established.

We now apply the localization by holomorphic two-form to this moduli
space. First, the form $\theta_+$ allows us to represent the virtual
cycle $\delta$ of $\fM_{0,2}(U^{[\ell]},d)$ as a homology class in
the Borel-Moore homology group $H\lsta^{BM}(\Lambda_{\theta_+})$, of
the set $\Lambda_{\theta_+}$ that is defined in (\ref{3.2}). By the
explicit construction of $\Lambda_{\theta_+}$, the Cartesian product
$$\begin{CD}
\Lambda_{\theta_+}\times_{U^{[\ell]}}(\fa_{-\ell}(E)\vac) @>>>
\Lambda_{\theta_+}
\\
@VVV @V{ev_1}VV\\
\fa_{-\ell}(E)\vac @>{\iota}>> U^{[\ell]}\\
\end{CD}
$$
is compact. Here $\iota$ is the tautological embedding and $ev_1$ is
the morphism defined by evaluating on the first marked point. Hence
because $U^{[\ell]}$ is smooth, the Gysin map
$$\iota\sta: H\lsta^{BM}(\Lambda_{\theta_+})\lra H\lsta\bl
\Lambda_{\theta_+}\times_{U^{[\ell]}}(\fa_{-\ell}(E)\vac)\br
$$
sends a Borel-Moore homology class to ordinary homology
class.

We let $\jmath\sta$ be the Gysin map associated to the taugological
inclusion and the second evaluation map:
$$\begin{CD}
\Lambda_{\theta_+}\times_{U^{[\ell]}}(\fa_{-\ell}(E)\vac) @>>>
\Lambda_{\theta_+}\times_{U^{[\ell]}}(\fa_{-\ell}(E)\vac)
\\
@VVV @V{ev_2}VV\\
(\fa_{-\ell}(X)\vac) @>{\jmath}>> U^{[\ell]}.\\
\end{CD}
$$
Note that $(\fa_{-\ell}(E)\vac) $ is a closed subset of $(\fa_{-\ell}(X)\vac)$, and thus
\begin{eqnarray*}
\Lambda_{\theta_+}\times_{U^{[\ell]}}(\fa_{-\ell}(E)\vac)\times_{U^{[\ell]}}(\fa_{-\ell}(X)\vac)=\Lambda_{\theta_+}\times_{U^{[\ell]}}(\fa_{-\ell}(E)\vac).
\end{eqnarray*}
Then an easy dimension count gives us
$$\jmath\sta\iota\sta \delta\in H_0(
\Lambda_{\theta_+}\times_{U^{[\ell]}}(\fa_{-\ell}(E)\vac)).
$$
We let $c_{\ell,d}$ be the degree of this cycle.

We remark that by the construction of the localized virtual cycle
$\delta$ and by the property of the Gysin maps, the so defined
number is universal in the sense that it does not depend on the
choice of the surface $U$, the form $\theta_+$ and the curve $C$, so
long as $\theta_+\upmo(0)$ intersects $E$ transversally at a single
point.

We define another universal constant $c_{\ell,d}^-$ similarly. We
keep the surface $U$, the curve $C$, but replace the holomorphic two
form $\theta_+$ by a meromorphic two-form $\theta_-$ that has no
vanishing divisor and has a smooth pole divisor $E$ that intersects
$C$ transversally at a single point. Then we take localized virtual
cycle $\delta_-\in H\lsta^{BM}(\Lambda_{\theta_-})$ of
$\fM_{0,2}(U^{[\ell]},d)$, and define $c_{\ell,d}^-$ to be the degree
of the class
$$\jmath\sta\iota\sta \delta_-\in H_0(
\Lambda_{\theta_-}\times_{U^{[\ell]}}(\fa_{-\ell}(E)\vac)).
$$

To proceed, we represent the Poincar\'e dual of $\be$ as
$$P.D\upmo(\be)=\alpha_1[E_1]-\alpha_2[E_2]+\alpha_0[E_0],
$$
where $E_1$ and $E_2$ are two smooth very ample divisors, $E_0\sub
X$ is a Riemann surface disjoint from $D_0\cup D_\infty$, and
$\alpha_1, \alpha_2$ are non-negative rational numbers. Then because
$$\blangle \mathfrak a_{-\ell}([E_0])\vac, \mathfrak
a_{-\ell}([X])\vac\brangle_{0, 2, d}=0,
$$
by the linearity of GW-invariants, we have
$$
\blangle \mathfrak a_{-\ell}(\be)\vac, \mathfrak
a_{-\ell}([X])\vac\brangle_{0,2,  d}= \sum_{i=1}^2
(-1)^{i-1}\alpha_i\blangle \mathfrak a_{-\ell}([E_i])\vac, \mathfrak
a_{-\ell}([X])\vac\brangle_{0, 2, d}.
$$
Then since we can arrange $E_i$ to intersects $D_0$ and $D_\infty$
transversally, the above sum is
$$\bl(\alpha_1[E_1]-\alpha_2[E_2])\cdot[D_0]\br c_{\ell,d} +
\bl(\alpha_1[E_1]-\alpha_2[E_2])\cdot[D_\infty]\br c_{\ell,d}^-.
$$
In case $c_{\ell,d}^-=-c_{\ell,d}$, then it becomes
$$\bl(\alpha_1[E_1]-\alpha_2[E_2])\cdot([D_0]-[D_\infty])\br c_{\ell,d}
=(\be\cdot c_1(K_X)) c_{\ell,d},
$$
as desired.

The identity $c_{\ell,d}^-=-c_{\ell,d}$ is easy to see. We let $X$
be a smooth K3 surface. Since its Hilbert scheme is holomorphic
symplectic, all its GW-invariants vanish. On the other hand, since
$K_X$ is trivial, we can find a meromorphic two form $\theta$ so
that $D_0$ and $D_\infty$ are non-empty. Following what we just
proved, say take $E_2=\emptyset$, we have
$$0=[E_1]\cdot[D_0] c_{\ell,d}+[E_1]\cdot [D_\infty] c_{\ell,d}^-
$$
for all smooth divisor $E_1$. This proves $c_{\ell,d}^-=c_{\ell,d}$,
and thus the Lemma.
\end{proof}

From the proof of the previous Proposition and the Lemma above, we
can get the following result.
\begin{corollary}\label{uni-cor} Let $\lambda$ and $\mu$ be as in
Corollary \ref{4.2}, then
\begin{eqnarray*}
&&\sum_{d\ge 0}\blangle A^\lambda_{\be}, c_1(\mathcal
O^{[n]}),A^\mu_{\bc}
\brangle_{0,3,d}\, q^d\\
=&&A^{\lambda}_{\be}\cup c_1(\mathcal O^{[n]})\cup A^{\mu}_{\bc}
+\sum_{\mu^3_i=\lambda^2_j=\ell}\blangle
A^{\lambda-\lambda^2_j}_{\be}, A^{\mu-\mu^3_i}_{\bc}\brangle
\sum_{d>0}d\,c_{\ell, d}\blangle \be_j, K_X\brangle\, q^d.
\end{eqnarray*}
\end{corollary}

\section{The quantum first Chern class operator}
\label{qfchern}

The Hilbert scheme $X^{[n]}$ admits a universal subscheme
\begin{eqnarray*}
Z_n\subset X^{[n]}\times X,\quad Z_n=\{(\xi, x)\,|\, x\in
\Supp(\xi)\}. \end{eqnarray*}

The sheaf $\pi_{1*}(\mathcal O_{Z_n})$ on $X^{[n]}$, where $\pi_1$
is the first factor projection, is a locally free sheaf of rank $n$.
We use $\mathcal O^{[n]}$ to denote this sheaf. The first Chern
class $c_1(\mathcal O^{[n]})$, treated as an operator on the
cohomology ring $H^*(X^{[n]})$ via the cup product, plays the
fundamental role in determining the ring structure of $H^*(X^{[n]})$
(see \cite{Lehn, LQW, C-G}).  Naturally, the action of $c_1(\mathcal
O^{[n]})$ on $H^*(X^{[n]})$ via the quantum cup product should play
the equally important role in the quantum cohomology ring
$H^*_\pi(X^{[n]})$. We call $c_1(\mathcal O^{[n]})\cup_\pi$ the quantum
first Chern class operator on $H^*_\pi(X^{[n]})$ when it acts via
the quantum product defined in subsection \S\ref{subsection-GW}.

When $X=\mathbb C^2$, the equivariant quantum first Chern class is
determined in \cite{O-P} via localization technique. The main task
of the rest of the paper is to determine the quantum first Chern
class operator for simply-connected surfaces.

 Consider the operator

\begin{eqnarray*}\label{M(q)}
M(q)&=&\sum_{k>0}\left(k\frac{(-q)^k}{(-q)^k-1}-\frac{q}{1+q}\right)
\mathfrak a_{-k}\mathfrak a_{k}(\tau_{2*}[K_X])\\
&&-\sum_{k>0}\frac{k-1}{2}\mathfrak a_{-k}\mathfrak
a_{k}(\tau_{2*}[K_X])-  \frac{1}{2} \sum_{k, \ell>0}\big(\mathfrak
a_{-k}\mathfrak a_{-\ell}\mathfrak a_{k+\ell}+\mathfrak
a_{-k-\ell}\mathfrak a_k\mathfrak a_{\ell}\big)(\tau_{3*}[X]).
\end{eqnarray*}
Here for $k\ge 1$, $\tau_{k\lsta}\colon H^*(X)\to H^*(X^k)$ is the linear
map induced by the diagonal embedding $\tau_k\colon X\to X^k$, and $\mathfrak a_{m_1}\ldots \mathfrak a_{m_k}(\tau_{k\lsta}(\alpha))$ denotes $\sum\limits_j\mathfrak a_{m_1}(\alpha_{j, 1})\ldots \mathfrak a_{m_k}(\alpha_{j, k})$ when $\tau_{k\lsta}\alpha=\sum\limits_j\alpha_{j, 1}\otimes \ldots\otimes\alpha_{j. k}$ via the K\"unneth decomposition of $H^*(X^k)$.

From Lehn's result \cite{Lehn} (see \cite{Q-W} also), $M(0)$ is the
first Chern class operator,
\begin{eqnarray*}
M(0)(A^\mu_{\bc})=c_1(\mathcal O_X^{[n]})\cup A^\mu_{\bc}.
\end{eqnarray*}

To calculate $M(q)(A^\mu_{\bc})$, it suffices to carry out the
following computations:
\begin{eqnarray}\label{three-terms}
\mathfrak a_{-k}\mathfrak a_k(\tau_{2*}[K_X])(A^\mu_{\bc}), \quad
a_{-k}\mathfrak a_{-\ell}\mathfrak
a_{k+\ell}(\tau_{3*}[X])(A^\mu_{\bc}), \quad  \mathfrak
a_{-k-\ell}\mathfrak a_{k}\mathfrak
a_{\ell}(\tau_{3*}[X])(A^\mu_{\bc}).
\end{eqnarray}

For the first term in (\ref{three-terms}), we have

\begin{eqnarray}
&&\mathfrak a_{-k}\mathfrak a_{k}(\tau_{2*}[K_X])A^\mu_{\bc}\nonumber\\
 =&&\sum_{\scriptstyle i=1\atop\scriptstyle \mu^2_i=k}^s(-\mu_i^2)\mathfrak
 a_{-k}([K_X\cdot C_i])A^{\mu-\mu^2_i}_{\bc}  \label{M"expan3}\\
 &&+\sum_{\scriptstyle i=1\atop\scriptstyle \mu^3_i=k}^t(-\mu_i^3)\mathfrak
 a_{-k}([K_X])A^{\mu-\mu^3_i}_{\bc}  \label{M"expan4}
 \end{eqnarray}

For the second term in (\ref{three-terms}), we have
\begin{eqnarray*}\label{M"expan1}
&&\mathfrak a_{-k}\mathfrak a_{-\ell}\mathfrak
a_{k+\ell}(\tau_{3*}[X])A^{\mu}_\bc\nonumber\\
 =&&\sum_{\scriptstyle i=1 \atop \scriptstyle k+\ell=\mu_i^1}^r(-\mu^1_i)\mathfrak a_{-k}\mathfrak
 a_{-\ell}(\tau_{2*}[q_i])A^{\mu-\mu^1_i}_{\bc}
 +\sum_{\scriptstyle i=1\atop\scriptstyle k+\ell=\mu_i^2}^s(-\mu^2_i)\mathfrak a_{-k}\mathfrak
 a_{-\ell}(\tau_{2*}[C_i])A^{\mu-\mu^2_i}_{\bc}\nonumber\\
 &&+\sum_{\scriptstyle i=1\atop\scriptstyle k+\ell=\mu_i^3}^t(-\mu^3_i)\mathfrak a_{-k}\mathfrak
 a_{-\ell}(\tau_{2*}[X])A^{\mu-\mu^3_i}_{\bc}.
 \end{eqnarray*}

For the third term in (\ref{three-terms}), we have
\begin{eqnarray*}
&&\mathfrak a_{-k-\ell}\mathfrak a_k\mathfrak
a_{\ell}(\tau_{3*}[X])A^\mu_{\bc}\label{M"expan2}
\\
=&&\sum_{\scriptstyle i=1\atop\scriptstyle
\mu^1_i=\ell}^r(-\mu_i^1)\mathfrak
 a_{-k-\ell}\mathfrak a_k(\tau_{2*}[q_i])A^{\mu-\mu_i^1}_{\bc}
 +\sum_{\scriptstyle i=1\atop\scriptstyle
\mu^2_i=\ell}^s(-\mu_i^2)\mathfrak
 a_{-k-\ell}\mathfrak a_k(\tau_{2*}[C_i])A^{\mu-\mu_i^2}_{\bc}\nonumber\\
 &&+\sum_{\scriptstyle i=1\atop\scriptstyle \mu^3_i=\ell}^t(-\mu_i^3)\mathfrak
 a_{-k-\ell}\mathfrak a_k(\tau_{2*}[X])A^{\mu-\mu_i^3}_{\bc}\nonumber\\
 =&&\sum_{\scriptstyle i=1\atop\scriptstyle \mu^1_i=\ell}^r\sum_{\scriptstyle j=1\atop\scriptstyle
  \mu^3_j=k}^tk\ell\mathfrak
 a_{-k-\ell}([q_i])A^{\mu-\mu_i^1-\mu_j^3}_{\bc}
+\sum_{\scriptstyle i=1\atop\scriptstyle \mu^2_i=\ell, i\neq j}^s\sum_{\scriptstyle j=1\atop\scriptstyle
  \mu^2_j=k}^sk\ell\mathfrak
 a_{-k-\ell}([C_i\cdot C_j])A^{\mu-\mu_i^2-\mu_j^2}_{\bc}\nonumber\\
 &&+\sum_{\scriptstyle i=1\atop\scriptstyle \mu^2_i=\ell}^s\sum_{\scriptstyle j=1\atop\scriptstyle
  \mu^3_j=k}^tk\ell\mathfrak
 a_{-k-\ell}([C_i])A^{\mu-\mu_i^2-\mu_j^3}_{\bc}+\sum_{\scriptstyle i=1\atop\scriptstyle \mu^3_i=\ell}^t\sum_{\scriptstyle j=1\atop\scriptstyle
  \mu^1_j=k}^rk\ell\mathfrak
 a_{-k-\ell}([q_j])A^{\mu-\mu_j^1-\mu_i^3}_{\bc}\nonumber\\
 &&+\sum_{\scriptstyle i=1\atop\scriptstyle \mu^3_i=\ell}^t\sum_{\scriptstyle j=1\atop\scriptstyle
  \mu^2_j=k}^sk\ell\mathfrak
 a_{-k-\ell}([C_j])A^{\mu-\mu_j^2-\mu_i^3}_{\bc}+\sum_{\scriptstyle i=1\atop\scriptstyle \mu^3_i=\ell,i\neq j}^t\sum_{\scriptstyle j=1\atop\scriptstyle
  \mu^3_j=k}^tk\ell\mathfrak
 a_{-k-\ell}([X])A^{\mu-\mu_j^1-\mu_i^3}_{\bc}\nonumber
 \end{eqnarray*}

The following result is the analogue of the Proposition
\ref{reduction}.
\begin{proposition}\label{operator-prop}
If $\blangle A^\lambda_{\be}, M(q)(A^\mu_{\bc})\brangle$ is not a
constant function of $q$, then either
$\ell(\mu^3)=\ell(\lambda^1)+1$ and $\ell(\lambda^3)=\ell(\mu^1)$,
or $\ell(\mu^3)=\ell(\lambda^1)$ and
$\ell(\lambda^3)=\ell(\mu^1)+1$.

In addition, assume $\ell(\mu^3)=\ell(\lambda^1)+1$ and
$\ell(\lambda^3)=\ell(\mu^1)$. Then  $\lambda^3=\mu^1$ as
partitions, and there exists an integer $\ell=\mu^3_i=\lambda^2_j$
for some integers $i$ and $j$ such that the partition $\lambda^1$ is
obtained from $\mu^3$ with $\ell$ deleted, and the partition $\mu^2$
is obtained from $\lambda^2$ with $\ell$ deleted.

\end{proposition}
\begin{proof}
 One can check easily, as cohomology class, $M(q)(A^\mu_{\bc})$ is of
 cohomology degree $\deg A^\mu_{\bc}+2$. Now take $A^\lambda_{\be}$ with
 $\deg A^\lambda_{\be}=4n-2-\deg A^\mu_{\bc}$. Therefore
 \begin{eqnarray*}
 \ell(\lambda^1)-\ell(\lambda^3)+\ell(\mu^1)-\ell(\mu^3)+1=0,
 \quad\hbox{i.e.}\quad
 \big(\ell(\lambda^3)-\ell(\mu^1)\big)+\big(\ell(\mu^3)-\ell(\lambda^1)\big)=1.
 \end{eqnarray*}

Since the second term and the third term in (\ref{three-terms}) only
contribute to the constant term of $\blangle A^\lambda_{\be},
M(q)(A^\mu_{\bc})\brangle$, we only need to consider the first term
in (\ref{three-terms}).

If $\ell(\lambda^3)-\ell(\mu^1)<0$, then $\ell(\mu^3)\ge
 \ell(\lambda^1)+2$.
 Then, for the terms in the summation of (\ref{M"expan3}) and (\ref{M"expan4}),  the number of terms $\mathfrak a_{-k}([X])$
 is
 more than $\ell(\lambda^1)$, thus $\blangle A^\lambda_{\be},
 M(q)(A^\mu_{\bc})\brangle=0$.

Similarly, by symmetry, if $\ell(\mu^3)-\ell(\lambda^1)<0$, we also
have $\blangle A^\lambda_{\be},
 M(q)(A^\mu_{\bc})\brangle=0$.

 Therefore if $\blangle A^\lambda_{\be},
 M(q)(A^\mu_{\bc})\brangle$ is not a constant function of $q$, we must have $\ell(\lambda^3)-\ell(\mu^1)\ge 0$ and
 $\ell(\mu^3)-\ell(\lambda^1)\ge 0$. Thus either
$\ell(\mu^3)=\ell(\lambda^1)+1$ and $\ell(\lambda^3)=\ell(\mu^1)$,
or $\ell(\mu^3)=\ell(\lambda^1)$ and
$\ell(\lambda^3)=\ell(\mu^1)+1$.

Assume
 $\ell(\mu^3)-\ell(\lambda^1)=1$ and $\ell(\lambda^3)=\ell(\mu^1)$.
 If $\blangle A^\lambda_{\be},
 M(q)(A^\mu_{\bc})\brangle$ is not a constant function of $q$, then $\blangle A^\lambda_{\be},
 M^{\prime}(q)(A^\mu_{\bc})\brangle$ cannot be a constant function of $q$, where
 \begin{eqnarray*}
 M^\prime(q)=M(q)-M(0).
 \end{eqnarray*}

The contribution from the terms in the summation of   (\ref{M"expan3}) must be zero since the number of terms
$\mathfrak a_{-k}([pt])$ there is $\ell(\mu^1)+1>\ell(\lambda^3)$.

Each  term in the summation of  (\ref{M"expan4})  has one more $\mathfrak
a_{-k}([\text{curve}])$ term and one less $\mathfrak a_{-\ell}([X])$ term
than $A^\mu_{\bc}$.
Thus by the Heisenberg commutation relation, we prove the
Proposition.
\end{proof}

\begin{corollary}\label{cor-qfc}
\begin{eqnarray*}
&&\blangle A^\lambda_{\be}, M(q)(A^\mu_{\bc})\brangle-\blangle A^\lambda_{\be}, M(0)(A^\mu_{\bc})\brangle \\
=&&\sum_{\scriptstyle1\le i\le t,\atop{\scriptstyle 1\le j\le
b,\atop \mu^3_i=\lambda^2_j=\ell}}\blangle
A^{\lambda-\lambda^2_j}_{\be}, A^{\mu-\mu^3_i}_{\bc}\brangle\blangle
E_j, K_X\brangle
(-1)^\ell\ell^2\left(\frac{\ell(-q)^\ell}{(-q)^\ell-1}-\frac{q}{1+q}\right).
\end{eqnarray*}
\end{corollary}

\section{\bf The projective plane}
\label{plane}

In Proposition \ref{reduction} in \S \ref{universality}, we see that
the two point extremal Gromov-Witten invariants can be reduced to
the computation of the case $\blangle \mathfrak
a_{-n}([C])\vac,\mathfrak a_{-n}([X])\vac \brangle_{0,2, d}$ where
$C$ is a curve.
 By Lemma \ref{uni-constant}, it suffices to carry out the
 computation for a particular surface. We choose the projective
 plane.  Since it is a toric surface, we can use computations in \cite{O-P}
 for the affine plane.

The main purpose of this section is to prove the following result.

\begin{proposition}\label{plane-prop}
Let $L$ be a line on the projective plane $X$. We have
\begin{eqnarray}\label{plane-prop-f1}
&&\sum_{d=0}^\infty\blangle \mathfrak a_{-n}([L])\vac, c_1(\mathcal
O^{[n]}), \mathfrak a_{-n}([X])\vac\brangle_{0,3,  d}q^d \nonumber\\
=&&\blangle\mathfrak a_{-n}([L])\vac, M(q)(\mathfrak
a_{-n}([X])\vac)\brangle.
\end{eqnarray}
\end{proposition}

Since $X$ is a toric surface, we can use the localization technique
to compute the corresponding equivariant Gromov-Witten invariants.
Since both sides of (\ref{plane-prop-f1}) are non-equivariant, we
get the conclusion of the Proposition \ref{plane-prop}. The
equivariant set-up is as follows.

 Let $[z_1, z_2, z_2]$ represent a point in the
projective space $X=\mathbb P^2$, and  $\mathbb T=\mathbb
C^*\times\mathbb C^*\times \mathbb C^*$ act on $X$ by
\begin{eqnarray*}
(s_1, s_2, s_3)\cdot [z_1, z_2, z_3]=[s_1z_1, s_2z_2, s_3z_3],
\quad\hbox{for }(s_1, s_2, s_3)\in \mathbb T.
\end{eqnarray*}

We have $H^*_{\mathbb T}(pt)=\mathbb C[t_1, t_2, t_3]$.

Let $Y$ be a $\mathbb T$-stable subvariety of $X$. We use $[Y]$ to
represent the equivariant class $Y\times_{\mathbb T}E\mathbb T$. By
the abuse of notation, we also use $[Y]$ to represent the
corresponding dual cohomology class $D^{-1}[Y]$ where $D\colon
H^k_{\mathbb T}(X)\to H^{\mathbb T}_k(X)$ is the Poincar\'e duality
morphism. Note the unusual convention on the degree of the
equivariant homology: if $Y$ has real codimension $k$ in $X$, then
$[Y]$ is a class in $H^{\mathbb T}_k(X)$ (see \cite{Vas}).

 There are three $\mathbb T$-fixed points:
\begin{eqnarray*}
q_1=[1,0,0],\quad q_2=[0,1,0],\quad q_3=[0,0,1].
\end{eqnarray*}
At $q_1$, under the identification  $[z_1,
z_2,z_3]=[1,z_2/z_1,z_3/z_1]$, the group $\mathbb T$ acts as
\begin{eqnarray*}
(s_1,s_2, s_3)\cdot (z_2/z_1,z_3/z_1)=(s_2s_1^{-1}z_2/z_1,
s_3s_1^{-1}z_3/z_1)\quad \hbox{for } (s_1, s_2, s_3)\in \mathbb T.
\end{eqnarray*}
The normal bundle $N_1$ of $q_1$ in $ X$, as a $\mathbb T$-module,
is isomorphic to $T_1^{-1}T_2\oplus T_1^{-1}T_3$, where $T_i$ is the
one-dimensional  representation given by $(s_1,s_2, s_3)\to s_i$.
Similarly, we also have
\begin{eqnarray*}
N_2=T_2^{-1}T_1\oplus T_2^{-1}T_3,\quad N_3=T_3^{-1}T_1\oplus
T_3^{-1}T_2
\end{eqnarray*}
as $\mathbb T$-modules,  where $N_i$ is the normal bundle of $q_i$
in $X$ regarded as a $\mathbb T$-module.

Let $L_1$ be the line in $X$ passing through $q_2$ and $q_3$, $L_2$
be the line passing through $q_1$ and $q_3$, and $L_3$ be the line
passing through $q_1$ and $q_2$. Near $q_1$, $L_3$ is given by the
equation $z_3/z_1=0$ and near $q_2$, $L_3$ is given by the equation
$z_3/z_2=0$. Thus by the localization, in the localized equivariant
cohomology
 \begin{eqnarray*}H^*_{\mathbb T}(X)^{\prime}=H^*_{\mathbb
T}(X)\otimes_{\mathbb C[t_1, t_2, t_3]}\mathbb C(t_1, t_2, t_3),
\end{eqnarray*}
we have
\begin{eqnarray}\label{L3}
[L_3]=\frac{(t_2-t_1)[q_1]}{(t_2-t_1)(t_3-t_1)}+\frac{(t_1-t_2)[q_2]}{(t_1-t_2)(t_3-t_2)}
=\frac{[q_1]}{(t_3-t_1)}+\frac{[q_2]}{(t_3-t_2)}.
\end{eqnarray}

Similarly we have
\begin{eqnarray}\label{X}
&&[L_2]=\frac{[q_1]}{t_2-t_1}+\frac{[q_3]}{t_2-t_3},\nonumber\\
 &&[L_1]
 =\frac{[q_2]}{t_1-t_2}+\frac{[q_3]}{t_1-t_3},\nonumber\\
 &&[X]
 =\frac{[q_1]}{(t_3-t_1)(t_2-t_1)}+\frac{[q_2]}{(t_3-t_2)(t_1-t_2)}+\frac{[q_3]}{(t_1-t_3)(t_2-t_3)}.
\end{eqnarray}

The anti-canonical class $-K_X$ can be written as, via localization,
\begin{eqnarray}\label{canonical}
&&[-K_X]\nonumber\\
&=&[L_1]+[L_2]+[L_3]\nonumber\\
&=&\frac{[q_2]}{t_1-t_2}+\frac{[q_3]}{t_1-t_3}+\frac{[q_1]}{t_2-t_1}+
\frac{[q_3]}{t_2-t_3}+\frac{[q_1]}{(t_3-t_1)}+\frac{[q_2]}{(t_3-t_2)}\nonumber\\
&=&\frac{(t_2+t_3-2t_1)[q_1]}{(t_2-t_1)(t_3-t_1)}+\frac{(t_1+t_3-2t_2)[q_2]}{(t_1-t_2)(t_3-t_2)}
+\frac{(t_1+t_2-2t_3)[q_3]}{(t_1-t_3)(t_2-t_3)}.
\end{eqnarray}
We also have, by the excess intersection formula
\begin{eqnarray*}
[q_1]\cdot [q_1]=(i_{1*}[q_1])\cdot
(i_{1*}[q_1])=i_{1*}\big([q_1]\cdot i_1^*i_{1*}[q_1]\big)=e^{\mathbb
T}(N_1)[q_1]=(t_3-t_1)(t_2-t_1)[q_1].
\end{eqnarray*}
Here $i_k\colon q_k\to X$ is the embedding. We don't distinguish
between $[q_1]$ with $i_{1*}[q_1]$ unless it is necessary.

Let's introduce a convention. Let $Y_i$ be $\mathbb T$-stable
subvarieties of $X$, $1\le i\le k$. We write
\begin{eqnarray*}
\mathfrak a_{n_1}\ldots\mathfrak a_{n_k}([Y_1\times\ldots\times
Y_k])\colon=\mathfrak a_{n_1}([Y_1])\ldots\mathfrak a_{n_k}([Y_k]) .
\end{eqnarray*}

Recall the  operator
\begin{eqnarray}\label{firstchernclass}
M(q) & =&\sum_{k>0}\big(\frac{k}{2}
         \frac{(-q)^k+1}{(-q)^k-1}-\frac{1}{2}\frac{(-q)+1}{(-q)-1}\big)
         \mathfrak a_{-k}\mathfrak a_k(\tau_{2*}[K_X])\nonumber\\
     &&~~~~~~~~~~~~~~~~~~~~~~~~~~~~~~~~    -\frac{1}{2}\sum_{k,\ell>0}\big(\mathfrak a_{-k}
\mathfrak a_{-\ell}\mathfrak a_{k+\ell}+\mathfrak
a_{-k-\ell}\mathfrak a_k\mathfrak a_{\ell}\big)(\tau_{3*}[X]).
\end{eqnarray}

With all the notations ready,  let's prove the Proposition
\ref{plane-prop}.
\begin{proof}[Proof of Proposition \ref{plane-prop}.]  Take $L=L_3$.
 Given $d>0$, we have
 \begin{eqnarray}\label{*1}
 &&\blangle \mathfrak a_{-n}([L_3])\vac, c_1(\mathcal
 O^{[n]}),\mathfrak a_{-n}([X])\vac \brangle_{0,3,d}\nonumber\\
 =&&
 d\blangle \mathfrak a_{-n}([L_3])\vac, \mathfrak a_{-n}([X])\vac
 \brangle_{0,2,d}\nonumber\\
 =&&d\Big\langle\mathfrak a_{-n}
 \left(\frac{[q_1]}{(t_3-t_1)}+\frac{[q_2]}{(t_3-t_2)}\right)
 \vac,\nonumber\\
 &&~~~~~~~~~~~\mathfrak a_{-n}\left(\frac{[q_1]}{(t_3-t_1)(t_2-t_1)}+
 \frac{[q_2]}{(t_3-t_2)(t_1-t_2)}+\frac{[q_3]}{(t_1-t_3)(t_2-t_3)}\right)
 \vac\Big\rangle_{0,2,d}\nonumber\\
 =&&d\frac{1}{(t_3-t_1)^2(t_2-t_1)}\blangle\mathfrak a_{-n}([q_1])\vac,\mathfrak a_{-n}([q_1])\vac\brangle_{0, 2, d}\nonumber\\
 &&+d\frac{1}{(t_3-t_2)^2(t_1-t_2)}\blangle\mathfrak a_{-n}([q_2])\vac,\mathfrak a_{-n}([q_2])\vac\brangle_{0,2,d}.
\end{eqnarray}

Let's introduction following notations for convenience.
\begin{eqnarray*}
M_1^\prime(q)&=&\frac{t_2+t_3-2t_1}{(t_2-t_1)(t_3-t_1)}\sum_{k>0}\big(\frac{k}{2}
         \frac{(-q)^k+1}{(-q)^k-1}-\frac{1}{2}\frac{(-q)+1}{(-q)-1}\big)\mathfrak a_{-k}([q_1])\mathfrak a_k([q_1]),\\
M_2^\prime(q)&=&\frac{t_1+t_3-2t_2}{(t_1-t_2)(t_3-t_2)}\sum_{k>0}\big(\frac{k}{2}
         \frac{(-q)^k+1}{(-q)^k-1}-\frac{1}{2}\frac{(-q)+1}{(-q)-1}\big)\mathfrak a_{-k}([q_2])\mathfrak a_k([q_2]),\nonumber\\
M_3^\prime(q)
&=&\frac{t_1+t_2-2t_3}{(t_1-t_3)(t_2-t_3)}\sum_{k>0}\big(\frac{k}{2}
         \frac{(-q)^k+1}{(-q)^k-1}-\frac{1}{2}\frac{(-q)+1}{(-q)-1}\big)\mathfrak a_{-k}([q_3])\mathfrak
         a_k([q_3]).\\
M_1^{\prime\prime}
&=&\frac{1}{2}\sum_{k,\ell>0}\frac{1}{(t_3-t_1)(t_2-t_1)}\big(\mathfrak
a_{-k}([q_1]) \mathfrak a_{-\ell}([q_1])\mathfrak
a_{k+\ell}([q_1])-\mathfrak a_{-k-\ell}([q_1])\mathfrak a_k([q_1])
\mathfrak a_{\ell}([q_1])\big),\\
M_2^{\prime\prime}
&=&\frac{1}{2}\sum_{k,\ell>0}\frac{1}{(t_3-t_2)(t_1-t_2)}\big(\mathfrak
a_{-k}([q_2]) \mathfrak a_{-\ell}([q_2])\mathfrak
a_{k+\ell}([q_2])-\mathfrak a_{-k-\ell}([q_2])\mathfrak a_k([q_2])
\mathfrak a_{\ell}([q_2])\big),\\
M_3^{\prime\prime}
&=&\frac{1}{2}\sum_{k,\ell>0}\frac{1}{(t_2-t_3)(t_1-t_3)}\big(\mathfrak
a_{-k}([q_3]) \mathfrak a_{-\ell}([q_3])\mathfrak
a_{k+\ell}([q_3])-\mathfrak a_{-k-\ell}([q_3])
\mathfrak a_k([q_3])\mathfrak a_{\ell}([q_3])\big).\\
\end{eqnarray*}

The terms $\blangle\mathfrak a_{-n}([q_i])\vac,\mathfrak
a_{-n}([q_i])\vac\brangle_{0, 2, d}$ in (\ref{*1})  was calculated
in \cite{O-P}. In fact
\begin{eqnarray}\label{*2}
&&\sum_{d=1}^\infty d\blangle\mathfrak a_{-n}([q_i])\vac,\mathfrak
a_{-n}([q_i])\vac\brangle_{0, 2, d}q^d\nonumber\\
=&&\sum_{d=1}^\infty\blangle\mathfrak a_{-n}([q_i])\vac,c_1(\mathcal
O^{[n]}), \mathfrak
a_{-n}([q_i])\vac\brangle_{0, 3, d}q^d\nonumber\\
=&&\blangle\mathfrak
a_{-n}([q_i])\vac,(M^{\prime}_i(q)+M^{\prime\prime}_i-M^{\prime}_i(0))
\mathfrak a_{-n}([q_i])\vac\brangle.
\end{eqnarray}

Now we have

\begin{eqnarray}\label{*3}
 &&\sum_{k>0}\left(\frac{k}{2}
         \frac{(-q)^k+1}{(-q)^k-1}-\frac{1}{2}\frac{(-q)+1}{(-q)-1}\right)
         \mathfrak a_{-k}\mathfrak a_k(\tau_{2*}[-K_X])\nonumber\\
&=&\sum_{k>0}\left(\frac{k}{2}
         \frac{(-q)^k+1}{(-q)^k-1}-\frac{1}{2}\frac{(-q)+1}{(-q)-1}\right)\nonumber\\
        && \mathfrak a_{-k}\mathfrak a_k\left(\frac{t_2+t_3-2t_1}{(t_2-t_1)(t_3-t_1)}[q_1 \times q_1]
         +\frac{t_1+t_3-2t_2}{(t_1-t_2)(t_3-t_2)}[q_2 \times q_2]
         +\frac{t_1+t_2-2t_3}{(t_1-t_3)(t_2-t_3)}[q_3 \times q_3]\right)\nonumber\\
         &=&\sum_{k>0}\left(\frac{k}{2}
         \frac{(-q)^k+1}{(-q)^k-1}-\frac{1}{2}\frac{(-q)+1}{(-q)-1}\right)\nonumber\\
         &&\cdot
         \Big(\frac{t_2+t_3-2t_1}{(t_2-t_1)(t_3-t_1)}\mathfrak a_{-k}([q_1])\mathfrak a_k([q_1])
         +\frac{t_1+t_3-2t_2}{(t_1-t_2)(t_3-t_2)}\mathfrak a_{-k}([q_2])\mathfrak a_k([q_2])\nonumber\\
         &&~~~~
         +\frac{t_1+t_2-2t_3}{(t_1-t_3)(t_2-t_3)}\mathfrak a_{-k}([q_3])\mathfrak a_k([q_3])\Big)\nonumber\\
&=&M_1^\prime(q)+M_2^\prime(q)+M_3^\prime(q),
\end{eqnarray}

Similarly, we also have
\begin{eqnarray}\label{*4}
 &&\frac{1}{2}\sum_{k,\ell>0}\big(\mathfrak a_{-k}
\mathfrak a_{-\ell}\mathfrak a_{k+\ell}-\mathfrak a_{-k-\ell}\mathfrak a_k\mathfrak a_{\ell}\big)(\tau_{3*}[X])\nonumber\\
&=&\frac{1}{2}\sum_{k,\ell>0}\big(\mathfrak a_{-k} \mathfrak
a_{-\ell}\mathfrak a_{k+\ell}-\mathfrak a_{-k-\ell}\mathfrak
a_k\mathfrak a_{\ell}\big)(\frac{[q_1\times q_1\times
q_1]}{(t_3-t_1)(t_2-t_1)})\nonumber\\
&&+\frac{1}{2}\sum_{k,\ell>0}\big(\mathfrak a_{-k} \mathfrak
a_{-\ell}\mathfrak a_{k+\ell}-\mathfrak a_{-k-\ell}\mathfrak
a_k\mathfrak a_{\ell}\big)(\frac{[q_2\times q_2\times
q_2]}{(t_3-t_2)(t_1-t_2)})\nonumber\\
&&+\frac{1}{2}\sum_{k,\ell>0}\big(\mathfrak a_{-k} \mathfrak
a_{-\ell}\mathfrak a_{k+\ell}-\mathfrak a_{-k-\ell}\mathfrak
a_k\mathfrak a_{\ell}\big)(\frac{[q_3\times q_3\times
q_3]}{(t_2-t_3)(t_1-t_3)})\nonumber\\
&=&M_1^{\prime\prime}+M_2^{\prime\prime}+M_3^{\prime\prime}
\end{eqnarray}
The first equality comes from $\tau_{3*}(q_i)=q_i\times q_i\times
q_i$ and the localization formula for $X$ expressed in terms of
fixed points $q_i$.

Combination of the formulae (\ref{*1}), (\ref{*2}), (\ref{*3}) and
(\ref{*4}) gives the conclusion of the Proposition.
\end{proof}

\section{\bf General surfaces and applications}
\subsection{General surfaces.}
One-point extremal Gromov-Witten invariants on the Hilbert scheme
$X^{[n]}$ for a simply-connected projective surface $X$ are
computed in \cite{L-Q}. As a consequence, the $3$-point extremal
Gromov-Witten invariants are computed for $X^{[2]}$ and the Ruan's
Cohomological Crepant Resolution Conjecture holds in this case.

In this section, we will determine two-point extremal Gromov-Witten
invariants of $X^{[n]}$.

\begin{theorem}\label{theorem}
Let $X$ be a simply connected projective surface. Then
\begin{eqnarray}\label{formula}
\sum_{d\ge 0}\blangle A^\lambda_{\be}, c_1(\mathcal O^{[n]}),
A^\mu_{\bc}\brangle_{0, 3, d}q^d=\blangle A^\lambda_{\be}, M(q)
A^\mu_{\bc}\brangle.
\end{eqnarray}
\end{theorem}
\begin{proof}
Let $L(q)$  denote the left hand side of (\ref{formula}) and $R(q)$
denote the right hand side of (\ref{formula}).

If $\ell(\lambda^3)-\ell(\mu^1)+\ell(\mu^3)-\ell(\lambda^1)\neq 1$,
both $L(q)$ and $R(q)$ equal to zero for the cohomological degree
reason.

If $\ell(\lambda^3)-\ell(\mu^1)+\ell(\mu^3)-\ell(\lambda^1)= 1$, but
$\ell(\lambda^3)\neq \ell(\mu^1)$ and $\ell(\mu^3)\neq
\ell(\lambda^1)$, by Proposition \ref{reduction} and Proposition
\ref{operator-prop}, we have
\begin{eqnarray*}
L(q)=A^\lambda_{\be}\cup c_1(\mathcal O^{[n]})\cup
A^\mu_{\bc}=L(0),\quad R(q)=\blangle A^\lambda_{\be}, M(0)\cup
A^\mu_{\bc}\brangle=R(0).
\end{eqnarray*}

Now $L(0)=R(0)$ follows from Lehn's result in \cite{Lehn} (see
\cite{Q-W} as well).

Next let's assume without loss of generality that
$\ell(\mu^3)=\ell(\lambda^1)+1$ and $\ell(\lambda^3)=\ell(\mu^1)$.

If $\lambda^3\neq\mu^1$, by Proposition \ref{reduction} and
Proposition \ref{operator-prop}, both $L(q)$ and $R(q)$ are constant
functions of $q$. Therefore $L(q)=L(0)=R(0)=R(q)$.

If $\lambda^3=\mu^1$, from Corollary \ref{uni-cor} and Lemma
\ref{uni-constant}, we get the formula
\begin{eqnarray*}
L(q)=L(0)+\sum_{\scriptstyle1\le i\le t,\atop{\scriptstyle 1\le j\le
b,\atop \mu^3_i=\lambda^2_j=\ell}}\blangle
A^{\lambda-\lambda^2_j}_{\be}, A^{\mu-\mu^3_i}_{\bc}\brangle\blangle
E_j, K_X\brangle \left(\sum_{d>0} dc_{\ell, d}q^d\right).
\end{eqnarray*}

From the discussion in \S \ref{qfchern}, we get
\begin{eqnarray*}
R(q)&=& \blangle
A^\lambda_{\be}, M(0)A^\mu_{\bc}\brangle+\blangle A^\lambda_{\be}, \big(M(q)-M(0)\big)A^\mu_{\bc}\brangle\\
&=&\blangle A^\lambda_{\be}, M(0)A^\mu_{\bc}\brangle\\
&&~~~~~~+\sum_{\scriptstyle1\le i\le t,\atop{\scriptstyle 1\le j\le
b,\atop \mu^3_i=\lambda^2_j=\ell}}\blangle
A^{\lambda-\lambda^2_j}_{\be}, A^{\mu-\mu^3_i}_{\bc}\brangle\blangle
E_j, K_X\brangle
(-1)^\ell\ell^2\left(\frac{\ell(-q)^\ell}{(-q)^\ell-1}-\frac{q}{1+q}\right).
\end{eqnarray*}

Now we need to determine $c_{\ell, d}$ explicitly. Since $c_{\ell,
d}$ is independent of the surface $X$, it suffices to consider the
case $X=\mathbb P^2$, $A^\lambda_{\be}=\mathfrak a_{\ell}(L)$ and
$A^\mu_{\bc}=\mathfrak a_{\ell}(X)$ where $L$ is a line in $X$.  By
Proposition \ref{plane-prop}, we have
\begin{eqnarray*}
-3\sum_{d>0}dc_{\ell, d}q^d&=&\blangle \mathfrak a_{-\ell}(L)\vac,
\big(M(q)-M(0)\big)\mathfrak
a_{-\ell}([X])\vac\brangle\\
&=&-3(-1)^\ell\ell^2\left(\frac{\ell(-q)^\ell}{(-q)^\ell-1}-\frac{q}{1+q}\right).
\end{eqnarray*}

Combination of all the formulae above gives us $L(q)=R(q)$.
\end{proof}

Thus we see that the quantum first Chern class has an explicit
formula, i.e., it is the operator $M(q)$.

\subsection{Application to Ruan's conjecture.}

In \cite{Q-W}, a vertex algebraic study of the cohomology
$H^*_{CR}(X^{(n)})$ was carried out. There is an irreducible
Heisenberg action $\{\mathfrak p_i(\alpha)\}_{i\in\mathbb Z,
\alpha\in H^*(X)}$ on $H^*_{CR}(X^{(n)})$ with a highest weight
vector $\vac$. Therefore there is a natural isomorphism
\begin{eqnarray*}
\Phi\colon H^*(X^{[n]}) \longrightarrow H^*_{CR}(X^{(n)})
\end{eqnarray*}
 as vector spaces. There is a counterpart of the first
Chern class of the tautological bundle on $H^*_{CR}(X^{(n)})$
defined in \cite{Q-W}. This class $O^1(1_X, n)\in H^*_{CR}(X^{(n)})$
defines an operator $\mathfrak b$ on $H^*_{CR}(X^{(n)})$ via the
Chen-Ruan product, which plays the similar role as $c_1(\mathcal
O^{[n]})$. It has the expression
\begin{eqnarray*}
\mathfrak b=-  \frac{1}{2} \sum_{k, \ell>0}\big(\mathfrak
p_{-k}\mathfrak p_{-\ell}\mathfrak p_{k+\ell}+\mathfrak
p_{-k-\ell}\mathfrak p_k\mathfrak p_{\ell}\big)(\tau_{3*}[X]).
\end{eqnarray*}

Therefore, a consequence of the Conjecture \ref{Ruan} for the
operator $\mathfrak b$ on $H^*_{CR}(X^{(n)})$ and the quantum first
Chern class operator on $H^*_\pi(X^{[n]})$ is the following equation

\begin{eqnarray}\label{Quantum-CR}
\blangle A^\lambda_{\be}, c_1(\mathcal O^{[n]})\cup_\pi
A^\mu_{\bc}\brangle =\blangle \Phi(A^\lambda_{\be}), O^1(1_X,
n)\cup_{CR}\Phi(A^\mu_{\bc})\brangle.
\end{eqnarray}

Recall that  $c_1(\mathcal O^{[n]})\cup_\pi$ is the operator
$M(-1)$. Therefore the way to prove the formula (\ref{Quantum-CR}),
under the identification of Heisenberg operators $\mathfrak
a_k(\alpha)\to \mathfrak p_k(\alpha)$, is to prove the operator
\begin{eqnarray*}
M(-1)=-  \frac{1}{2} \sum_{k, \ell>0}\big(\mathfrak a_{-k}\mathfrak
a_{-\ell}\mathfrak a_{k+\ell}+\mathfrak a_{-k-\ell}\mathfrak
a_k\mathfrak a_{\ell}\big)(\tau_{3*}[X]).
\end{eqnarray*}

 Note that
$$\frac{k}{2}\frac{(-q)^k+1}{(-q)^k-1}-\frac{1}{2}\frac{(-q)+1}{(-q)-1}$$
is well defined at $q=-1$ by using  L'Hospital's rule. In fact,


\begin{eqnarray*}
\left(\frac{k}{2}\frac{(-q)^k+1}{(-q)^k-1}-\frac{1}{2}\frac{(-q)+1}{(-q)-1}\right)\big|_{q=-1}=0,
\end{eqnarray*}
 and therefore
$$M(-1)=-\frac{1}{2}\sum_{k, \ell>0}\big(\mathfrak
a_{-k}\mathfrak a_{-\ell}\mathfrak a_{k+\ell}+\mathfrak
a_{-k-\ell}\mathfrak a_k\mathfrak a_{\ell}\big)(\tau_{3*}[X]).$$

Thus we proved the formula (\ref{Quantum-CR}).

\begin{remark}
When the surface $X$ is $K3$, the combination of results in \cite{L-S} and \cite{F-G, Ur} verifies Ruan's conjecture.
When $X=\mathbb P^2$, it is is shown in \cite{ELQ}  that the three-point extremal GW-invariants for $X^{[3]}$ can be reduced to two-point extremal GW-invaraints. By the result in this paper, the three-point extremal GW-invariants for $X^{[3]}$ are determined.
\end{remark}

\end{document}